\documentclass[12pt,leqno]{article}
\usepackage{amsmath, amsthm, amsfonts, amssymb, color}
\usepackage{mathrsfs}
\theoremstyle{definition}

\newcommand{\scr}[1]{\mathscr #1}
\definecolor{wco}{rgb}{0.5,0.2,0.3}

\numberwithin{equation}{section} \theoremstyle{remark}

\newcommand{\ua}{\uparrow}

\title{
{\bf Degenerate SDEs in Hilbert Spaces with Rough Drifts}
\footnote{FW is supported in part by
NNSFC (11131003, 11431014), the 985 project and the Laboratory of Mathematical and  Complex Systems, XZ is supported partly by
NNSFC (11271294, 11325105).}}
\author{
{\bf Feng-Yu Wang$^{a),b)}$ and Xicheng Zhang$^{c)}$ }\\
\footnotesize{a) School of Mathematical Sciences,
Beijing Normal University, Beijing 100875, China}\\
 \footnotesize{b) Department of Mathematics,
Swansea University, Singleton Park, SA2 8PP, UK}\\
\footnotesize{c) School of Mathematics and Statistics,
Wuhan University, Wuhan 430072, China}}

\begin{document}
\def\tttext#1{{\normalfont\ttfamily#1}}
\def\R{\mathbb R}  \def\ff{\frac} \def\ss{\sqrt} \def\B{\mathbf B}
\def\N{\mathbb N} \def\kk{\kappa} \def\m{{\bf m}}
\def\dd{\delta} \def\DD{\Delta} \def\vv{\varepsilon} \def\rr{\rho}
\def\<{\langle} \def\>{\rangle} \def\ggm{\Gamma} \def\ggm{\gamma}
  \def\nn{\nabla} \def\pp{\partial} \def\EE{\scr E}
\def\d{\text{\rm{d}}} \def\bb{\beta} \def\aa{\alpha} \def\D{\scr D}
  \def\si{\sigma} \def\ess{\text{\rm{ess}}}
\def\beg{\begin} \def\beq{\begin{equation}}  \def\F{\scr F}
\def\Ric{\text{\rm{Ric}}} \def\Hess{\text{\rm{Hess}}}
\def\e{\text{\rm{e}}} \def\ua{\underline a} \def\OO{\Omega}  \def\oo{\omega}
 \def\tt{\tilde} \def\Ric{\text{\rm{Ric}}}
\def\cut{\text{\rm{cut}}} \def\P{\mathbb P}
\def\C{\scr C}     \def\E{\mathbb E}
\def\Z{\mathbb Z} \def\II{\mathbb I}
  \def\Q{\mathbb Q}  \def\LL{\Lambda}
  \def\B{\scr B}    \def\ll{\lambda}
\def\vp{\varphi}\def\H{\mathbb H}\def\ee{\mathbf e} \def\GG{\Gamma}\def\gg{\gamma}

\maketitle

\begin{abstract} The existence and uniqueness of  mild solutions are proved for   a class of degenerate stochastic
differential equations on Hilbert spaces where the drift is Dini continuous in  the component
with noise and H\"older continuous of order larger than $\ff 2 3$ in
the other component.
 In the finite-dimensional case the Dini continuity is further weakened.
The main results are applied to solve second order
stochastic  systems driven by   space-time white noises.

\end{abstract} \noindent


\section{Introduction}

Let $\H_i$ ($i=1,2,3$) be separable Hilbert spaces, and let $\scr L(\H_i;\H_j)$ be the class of all bounded linear operators from $\H_i$ to $\H_j \ (1\le i,j\le 3).$ We shall simply denote the norm and inner product by $|\cdot|$ and $\<\cdot,\cdot\>$ for Hilbert spaces, and let $\|\cdot\|$ stand for the operator norm.

Let $W_t$ be a cylindrical Brownian motion on $\H_3$; i.e. for an orthonormal basis $\{h_i\}_{i\ge 1}$ on $\H_3$, we have
$$W_t= \sum_{i\ge 1}B_t^i h _i,$$ where $\{B_t^i\}_{i\ge 1}$ is a family of independent one-dimensional Brownian motions.
Let $\{\F_t\}_{t\ge 0}$ be the natural filtration induced by $W_t$.

We consider the following degenerate stochastic evolution equation on $\H:=\H_1\times\H_2$:
\beq\label{1.1}\beg{cases} \d X_t= \big\{A_1X_t +BY_t\big\}\d t,\\
\d Y_t= \big\{b_t(X_t,Y_t)+A_2Y_t\big\}\d t +\si_t\d W_t,\end{cases}\end{equation}
where $ B\in \scr L(\H_2;\H_1),\ \si: [0,\infty)\to \scr L(\H_3;\H_2)$, $b: [0,\infty)\times\H\to \H_2$ are measurable and locally bounded,
and for every $i=1,2$, $(A_i,\D(A_i))$ is a bounded above linear operator  generating a strongly continuous   semigroup $\e^{tA_i}$ on $\H_i$. We will let $\nn, \nn^{(1)}$ and $\nn^{(2)}$ denote the gradient operators on $\H,\H_1$ and $\H_2$ respectively.

\beg{defn} A  continuous adapted process $(X_t,Y_t)_{t\in [0,\zeta)}$ is called a mild solution to \eqref{1.1} with life time $\zeta$, if $\zeta>0$ is
an $\F_t$-stopping time
such that $\P$-a.s. $\limsup_{t\uparrow \zeta} (|X_t|+|Y_t|)=\infty$ holds on $\{\zeta<\infty\}$ and, $\P$-a.s. for all $t\in [0,\zeta),$
\beg{equation*}\left\{\beg{aligned}
 & X_t= \e^{tA_1}X_0 +\int_0^t \e^{(t-s)A_1}BY_s \d s,\\
& Y_t=\e^{tA_2}Y_0 +\int_0^t \e^{(t-s)A_2} b_s(X_s,Y_s)\d s +\int_0^t \e^{(t-s)A_2}\si_s \d W_s.\end{aligned}\right.
\end{equation*}
\end{defn}

The purpose of this paper is to investigate the existence/uniqueness  of the mild solution
under some Dini's type continuity conditions on the drift $b$.
  The main idea is to construct a map which transforms the original equation into an equation with regular enough coefficients ensuring the pathwise uniqueness of  the solution.  This idea goes back to  \cite{V,Z} where finite-dimensional SDEs with singular drift are investigated, see also \cite{KR,Zh} for further developments. In recent years, this argument has been developed in \cite{DF,DR1,DR2,RN,W14} for non-degenerate SDEs in Hilbert spaces.  The main difficulty of the study for the present degenerate equation is that the semigroup $P_t^0$ associated to the linear equation (i.e. $b=0$) has worse gradient estimates with respect to $x\in \H_1$. More precisely, unlike in the non-degenerate case one has $\|\nn P_t^0\|_{\infty\to\infty}\approx t^{-1/2}$ for small $t>0$ which is integrable over $[0,1]$, for the present model one has $\|\nn^{(1)}P_t^0\|_{\infty\to\infty}
\approx t^{-3/2}$ which is not integrable, where $\|P\|_{\infty\to\infty}:=\sup_{\|f\|_\infty\le 1} \|Pf\|_\infty$ for a linear operator $P$ and the uniform norm $\|\cdot\|_\infty$.   To reduce the singularity for small $t>0$, we will use some other norms to replace  $\|\cdot\|_{\infty\to\infty}$. This leads to different type continuity conditions on $b$. Indeed, we will need the H\"older continuity of $b$ in the first component $x\in \H_1$, and a Dini type continuity of $b$ in the second component $y\in \H_2$ as in \cite{W14} where the non-degenerate equation is concerned.

To ensure the required gradient estimates on $P_t^0$, we make the following assumptions on the linear part.
\beg{enumerate} \item[{\bf (H1)}] $\si_t\si_t^*$ is invertible in $\H_2$ with locally bounded $\|(\si_t\si_t^*)^{-1}\|$ in $t\ge 0.$
\item[{\bf (H2)}] $BB^*$ is invertible in $\H_1$, and $B\e^{tA_2}= \e^{tA_1}\e^{tA_0}B$ for some   $A_0\in \scr L(\H_1,\H_1)$ and all $t\ge 0$.
\item[{\bf (H3)}] $-A_2$ is self-adjoint having  discrete spectrum $0<\ll_1\le\ll_2\le\cdots$ counting multiplicities
such that $\sum_{i\ge 1} \ff 1 {\ll_i^{1-\delta}}<\infty$ for some $\delta\in (0,1).$
  \end{enumerate}

Since $\si_t$ is locally bounded in $t\ge 0$, $B$ is bounded,  and $ A_1 $ is bounded above so that $\|\e^{A_1t}\|\le \e^{ct}$ holds for some constant $c\ge 0$, it is well known from \cite{DZ} that {\bf (H3)} implies the existence, uniqueness and non-explosion of a continuous mild solution to the linear equation,  i.e. \eqref{1.1}
with $b=0$. As It\^o's formula does not apply directly to the mild solution, in the study we will make finite-dimensional approximations.
Throughout the paper, we let $\{e_i\}_{i\ge 1}$ be the eigenbasis of $A_2$, which is an orthonormal basis in $\H_2$ such that $A_2e_i=-\ll_i e_i.$ For any $n\ge 1$, let
$\H_2^{(n)}=\text{span}\{e_1,\cdots, e_n\}$, and let $\pi_2^{(n)}: \H_2\to \H_2^{(n)}$ be the orthogonal projection. Next, let $\H_1^{(n)}= B\H_2^{(n)}$ and
$\pi_1^{(n)}:\H_1 \to \H_1^{(n)}$ be the orthogonal projection. Since $BB^*$ is invertible, we have $\lim_{n\to\infty} \pi_1^{(n)}x= x$ for $x\in \H_1$. Let
$$\pi^{(n)}=(\pi_1^{(n)}, \pi_2^{(n)}): \H\to \H_n:= \H_1^{(n)}\times\H_2^{(n)}.$$ In our study of finite-dimensional approximations, we will need the following   assumption which is trivial in the finite-dimensional setting.

\beg{enumerate} \item[{\bf (H4)}] There exists $n_0\ge 1$ such that for any $n\ge n_0$, $\pi_1^{(n)} B= B\pi_2^{(n)}$ on $\H_2$, and $\pi_1^{(n)}A_1 = A_1\pi_1^{(n)}$   on $\D(A_1)$.
  \end{enumerate}

We introduce the following classes of functions to characterize the continuity modulation of the drift $b$:
\beg{equation*}\beg{split}
&\D_0:= \Big\{\phi:  [0,\infty)\to [0,\infty)  \text{\ is\ increasing  with }\phi(0)=0 \text{ and } \phi(s)>0\ \text{for}\ s>0\Big\},\\
&\D_1:=\left\{\phi \in \D_0: \  \phi^2 \text{\ is\ concave and} \int_0^1\ff{\phi(s)}s\d s<\infty\right\}.\end{split}\end{equation*}
We remark that the condition $\int_0^1\ff{\phi(s)}s\d s<\infty$ is well known as
Dini's condition,
due to the notion of Dini's continuity.
Obviously, the class $\D_1$ contains $\phi(s):= \ff{K}{\{\log(c+s^{-1})\}^{1+r}}$ for constants $K,r>0$ and large enough $c\ge \e$ such that $\phi^2$ is concave.

\beg{thm}\label{T1.1} Assume  {\bf (H1)}-{\bf (H4)}.    \beg{enumerate}\item[$(1)$] If for any $n\ge 1$ there exist  $\aa_n\in (\ff 2 3,1),$
$\phi_n\in\D_1$ and a constant $K_n>0$ such that
$$|b_t(x,y)-b_t(x',y')|\le K_n |x-x'|^{\aa_n}+ \phi_n(|y-y'|),\ \ t\in [0,n], |(x,y)|\lor |(x',y')|\le n, $$
then for any  $(X_0,Y_0)\in\H$, the equation $\eqref{1.1}$ has a unique mild solution $(X_t,Y_t)_{t\in [0,\zeta)}$  with life time $\zeta$.
\item[$(2)$] If moreover
\beq\label{C1} \< b_t (x, y+y'), y\> \le \ell_t(|x|^2+|y|^2)+ h_t(|y'|),\ \ \ x\in\H_1, y,y'\in \H_2, t\ge 0 \end{equation}   holds for some increasing function $\ell,h: [0,\infty)\times [0,\infty)\to (0,\infty)$ such that $\int_1^\infty \ff {\d s}{\ell_t(s)}=\infty$ holds for all $t\ge 0$, then the unique mild solution is non-explosive, i.e. the life time $\zeta=\infty\ \P$-a.s.
\end{enumerate}\end{thm}

To illustrate this result, we consider the following example of second order stochastic   system driven by   white noise.

\paragraph{Example 1.1.} Let $D\subset \R^d$ be a bounded open domain, and let $\DD$ be the Dirichlet Laplace operator on $D$. Consider the  equation
\beq\label{We0}\beg{split}  \pp_t^2 u(t,x)= & b_t\big(u(t,x), \pp_t u(t,x)+ (-\DD)^\theta u(t,\cdot)(x)\big)
 -(-\DD)^{2\theta} u(t,\cdot)(x)\\
 & - 2 (-\DD)^\theta \pp_t u(t,\cdot)(x) +\ff{W(\d t,\d x)}{\d t\d x},\ \ t\ge 0, x\in D.\end{split}\end{equation}
Here, $\theta>0$ is a constant, $W$ is a Brownian sheet (space-time white noise) on $\R^d$, and $b: [0,\infty)\times \R^2\to \R$ is measurable such that for any $T>0$,
\beq\label{We1} |b_t(u,v)-b_t(u',v')|\le C\big(|u-u'|^\aa+ \phi(|v-v'|)\big),\ \ t\in [0,T], u,v,u',v'\in\R\end{equation} holds for some constants
$C>0, \aa\in (\ff 2 3,1)$, and some $\phi\in \D_1.$

To solve this equation using Theorem \ref{T1.1}, we take $\H_i= L^2(D;\d x)$ for $i=1,2,3$, and
$$W_t= \sum_{i=1}^\infty e_i\int_{[0,t]\times D} e_i(x)W(\d s,\d x)$$
for $\{e_i\}_{i\ge 1}$ the unitary eigenbasis of $\DD$. Letting
$$X_t= u(t,\cdot),\ \ Y_t= \pp_t u(t,\cdot)+ (-\DD)^\theta u(t,\cdot),$$ we   reformulate  \eqref{We0} as
\beq\label{1.1E}\beg{cases} \d X_t= \big\{Y_t- (-\DD)^\theta X_t \big\}\d t,\\
\d Y_t= \big\{b_t(X_t,Y_t)-(-\DD)^\theta Y_t\big\}\d t + \d W_t.\end{cases}\end{equation}
Obviously, assumptions {\bf (H1)}, {\bf (H2)} and {\bf (H4)} hold for $B=\si_t=I$ (the identity operator) and $A_1=A_2= -(-\DD)^\theta$. Moreover, since the eigenvalues of
$(-\DD)^\theta$ satisfy $\ll_i\ge ci^{2\theta/d}$ for some constant $c>0$ and all $i\ge 1$, assumption {\bf (H3)} holds for $A_2=-(-\DD)^\theta$ provided $\theta>\ff d 2.$ Finally, by Jensen's inequality it is easy to see that \eqref{We1} implies
$$\|b(f,g)-b(\tt f,\tt g)\|_{L^2(D)}\le C\big(\|f-\tt f\|_{L^2(D)}+ \phi(\|g-\tt g\|_{L^2(D)})\big),\ \ f,g,\tt f,\tt g\in L^2(D;\d x), t\in [0,T],$$ where the constant $C$ might be different  if the volume of $D$ is not equal to 1. Therefore, by Theorem \ref{T1.1}, for any $\theta>\ff d 2$ the equation \eqref{1.1E} has a unique mild solution on $L^2(D;\d x)\times L^2(D;\d x)$ which is non-explosive.

\

Next, we consider the finite-dimensional case, i.e. consider the following degenerate  SDE on $\R^m\times\R^d$:
\beq\label{1.2}\beg{cases} \d X_t= \big\{AX_t +BY_t\big\}\d t,\\
\d Y_t= b_t(X_t,Y_t)\d t +\si_t\d W_t,\end{cases}\end{equation}
 where $A$ is an $m\times m$-matrix, $B$ is an $m\times d$-matrix,   $\si: [0,\infty)\to \scr L(\H_3;\R^d)$ is measurable and locally bounded,   and $W_t$ is a cylindrical Brownian motion on $\H_3$. In this case, Theorem \ref{T1.1} can be improved by using the following larger class $\D_2$ to replace $\D_1$:
 $$\D_2:=\left\{\phi \in \D_0: \  \phi^2 \text{\ is\ concave,} \int_0^1\ff{\d t}{t\big(1+\int_t^1\ff{\phi(s)}s\d s\big)^2}=\infty\right\},$$ which includes $\phi(s):= \ff K {\ss{\log (c+ s^{-1})}}$ for some constant  $K>0$ and large enough constant $c>0.$

\beg{thm}\label{T1.2} Let $\H_1=\R^m$ and $\H_2=\R^d$ be finite-dimensional.  Assume that $BB^*$ and
 $\si_t\si_t^*$ are invertible with $(\si_t\si_t^*)^{-1}$ locally bounded in $t\ge 0$.  \beg{enumerate}\item[$(1)$] If for any $n\ge 1$ there exist  $\aa_n\in (\ff 2 3,1),$
$\phi_n\in\D_2$ and a constant $K_n>0$ such that
$$|b_t(x,y)-b_t(x',y')|\le K_n |x-x'|^{\aa_n}+ \phi_n(|y-y'|),\ \ t\in [0,n], |(x,y)|\lor |(x',y')|\le n, $$
then for any  $(X_0,Y_0)\in\R^{m+d}$, the equation $\eqref{1.2}$ has a unique   solution $(X_t,Y_t)_{t\in [0,\zeta)}$  with life time $\zeta$.
\item[$(2)$] If moreover
\beq\label{C1'} \< b_t (x, y), y\> \le \ell_t(|x|^2+|y|^2),\ \ \ x\in\R^m, y\in \R^d, t\ge 0 \end{equation}   holds for some increasing function $\ell,h: [0,\infty)\times [0,\infty)\to (0,\infty)$ such that $\int_1^\infty \ff {\d s}{\ell_t(s)}=\infty$ holds for all $t\ge 0$, then the unique mild solution is non-explosive, i.e. the life time $\zeta=\infty\ \P$-a.s. \end{enumerate}  \end{thm}

\paragraph{Example 1.2.}
Consider the following second order stochastic differential equation on $\R^{d}$:
$$
\ff{\d^2 X_t}{\d t^2}=b_t\Big(X_t,\ff{\d X_t}{\d t}\Big)+\sigma\dot W_t,
$$
where $ W_t$ is the $d$-dimensional
Brownian, 
$\sigma \in \R^d\otimes\R^d$ is  invertible,
and $b:[0,\infty)\times \R^d\times\R^d\to\R^d$ is measurable such that for any $T>0$ the condition \eqref{We1} holds   for some constants $C>0$, $\alpha\in(\tfrac{2}{3},1)$ and some function $\phi\in\D_2.$  By letting $m=d$ and $Y_t=\ff{\d X_t}{\d t}$, we reformulate this equation as
\beq\label{1.2E}\beg{cases} \d X_t=  Y_t\,\d t,\\
\d Y_t= b_t(X_t,Y_t)\d t +\si \d W_t.\end{cases}\end{equation}According to Theorem \ref{T1.2}, for any initial point this equation has   a unique solution which is non-explosive.

\

We would like to point out that in the finite-dimensional
setting, the pathwise uniqueness for equation \eqref{1.2}  with H\"older continuous drifts has been investigated in a preprint by Chaudru de Raynal
(http://hal.archives-ouvertes.fr/hal-00702532/document).
However, we found some gaps in the proof, for instance,
the probabilistic representation of the solution is wrongly used and this is crucial  in related calculations.
Obviously, in Theorem \ref{T1.2} the   condition on $b$ along the second component $y$ is much weaker than the H\"older continuity.

\

The remainder of the paper is organized as follows. In Section 2, we investigate gradient estimates on the  semigroup $P_{s,t}^0$ associated to the linear equation (i.e. $b=0$).
These gradient estimates are
then used in Section 3 to construct and study the regularization transform. In Section 4 we use the regularization transform to represent the mild solution to \eqref{1.1}, which enables us to prove Theorems \ref{T1.1} and Theorem \ref{T1.2} in Section 5.

\section{Gradient estimates on $P_{s,t}^0$}

For any $s\ge 0$, consider the linear equation
\beq\label{2.1}\beg{cases} \d X_{s,t}^0= \big\{A_1X_{s,t}^0 +BY_{s,t}^0\big\}\d t,\\
\d Y_{s,t}^0=  A_2Y_{s,t}^0\,\d t +\si_t\d W_t,\ \ t\ge s.\end{cases}\end{equation}
By {\bf (H1)}-{\bf (H3)} and Duhamel's formula,  the unique solution of this equation starting at $(x,y)\in \H$ at time $s$ is given by:
\beq\label{2.2} \left\{\beg{aligned}& X_{s,t}^0= \e^{(t-s)A_1} x +\int_s^t \e^{(t-r)A_1}BY_{s,r}^0\d r,\\
&Y_{s,t}^0= \e^{(t-s)A_2}y + \int_s^t\e^{(t-r)A_2}\si_r\d W_r.\end{aligned}\right.\end{equation} To indicate the dependence on the initial point, we also denote the solution by $(X_{s,t}^0, Y_{s,t}^0)(x,y).$ Let $P_{s,t}^0$ be the Markov operator
associated to $(X_{s,t}^0, Y_{s,t}^0)$ , i.e.
$$P_{s,t}^0 f(x,y)= \E f((X_{s,t}^0, Y_{s,t}^0)(x,y)),\ \ t\ge s\ge 0, (x,y)\in \H, f\in \B_b(\H).$$
By the Markov property, we have $P^0_{s,r}P_{r,t}^0= P_{s,t}^0$ for $0\le s\le r\le t.$ \

We first present a Bismut type derivative formula for $P_{s,t}^0$. Let
$$Q_t= \int_0^t s(t-s) \e^{sA_0}BB^*\e^{sA_0^*}\d s,\ \ t>0,$$ where $A_0$ is in {\bf (H2)}. Since $BB^*$ is invertible and $A_0$ is bounded, $Q_t^{-1}$ is invertible for every $t>0$, and for any $T>0$ there exists a constant $c>0$ such that
\beq\label{Q} \|Q_t^{-1}\|\le \ff{c} {t^3},\ \  \ t\in (0,T].\end{equation} Next,
for any  $s\in [0,T)$ and $v=(v_1,v_2)\in \H$, let
\beg{equation*}\beg{split} & V^v_{s,T} := Q_{T-s}^{-1}\bigg[v_1 +  \int_{s}^{T} \ff{T-r}{T-s}\e^{(r-s)A_0^*}Bv_2\,\d r \bigg],\\
& \Phi^v_{s,T} (r):= \e^{(r-s)A_2} \bigg[\ff{v_2}{T-s} +\ff{\d}{\d r}\big\{(r-s)(T-r)B^*\e^{(r-s)A_0^*}\big\}V^v_{s,T} \bigg].\end{split}\end{equation*}

\beg{thm}\label{T2.1} For any  $s\in [0,T), v=(v_1,v_2)\in \H$, $f\in \B_b(\H)$, and $(x,y)\in \H$,   there holds
\begin{align}\label{BS}
(\nn_v P_{s,T}^0 f)(x,y) = \E\bigg[ f\big((X_{s,T}^0, Y_{s,T}^0)(x,y)\big)\int_{s}^{T} \big\<\si_r^*(\si_r\si_r^*)^{-1} \Phi^v_{s,T}(r), \d W_r\big\>\bigg].
\end{align}
\end{thm}

\beg{proof} We use the argument of coupling by change of measures as in \cite{GW}
where the finite-dimensional case is considered. For any $\vv\in [0,1)$, let $(X_{s,t}^\vv, Y_{s,t}^\vv)_{t\ge s}$ solve the equation
\beq\label{2.3}\beg{cases}
\d X_{s,t}^\vv= \big\{A_1X_{s,t}^\vv +BY_{s,t}^\vv\big\}\d t,&X^\vv_{s,s}=x+\vv v_1,\\
\d Y_{s,t}^\vv=  \big\{A_2Y_{s,t}^\vv- \vv \Phi^v_{s,T}(t)\big\}\,\d t +\si_t\d W_t, &Y_{s,s}^\vv= y+\vv v_2.\end{cases}\end{equation}
Noticing that
\beg{equation*}\beg{cases} \d (X_{s,t}^\vv-X_{s,t}^0)= \big\{A_1(X_{s,t}^\vv-X_{s,t}^0)
+B(Y_{s,t}^\vv-Y_{s,t}^0)\big\}\d t,\\
\d (Y_{s,t}^\vv- Y_{s,t}^0)= \big\{A_2(Y_{s,t}^\vv- Y_{s,t}^0)- \vv \Phi^v_{s,T}(t)\big\} \d t,\end{cases}\end{equation*}
by Duhamel's formula and the definition of $\Phi^v_{s,T}$, we have
\beq\label{2.4} \beg{split}& Y_{s,t}^\vv- Y_{s,t}^0= \vv \e^{(t-s)A_2}v_2 -\vv \int_{s}^t \e^{(t-r)A_2} \Phi^v_{s,T}(r)\d r\\
&= \vv \e^{(t-s)A_2} \left[\ff{T-t}{T-s} v_2- (t-s)(T-t)B^*\e^{(t-s)A_0^*}V^v_{s,T}\right].
\end{split}\end{equation}
On the other hand, by {\bf (H2)}, we also have
\beq\label{2.5} \beg{split}& X_{s,t}^\vv-X_{s,t}^0  = \vv \e^{(t-s)A_1}v_1 +\int_{s}^t \e^{(t-r)A_1} B (Y_{s,r}^\vv- Y_{s,r}^0)\d r\\
& =\vv \e^{(t-s)A_1}\left[v_1 +\int_s^t\e^{(r-s)A_0}\Big(\ff{T-r}{T-s}Bv_2 -(r-s)(T-r)   BB^* \e^{(r-s)A_0^*} V^v_{s,T}\Big)\d r \right].\end{split}\end{equation}
In particular, by the definition of $V^v_{s,T}$, \eqref{2.4} and \eqref{2.5} imply
\beq\label{2.6}
(X_{s,T}^\vv, Y_{s,T}^\vv) = (X_{s,T}^0, Y_{s,T}^0),\ \ \vv \in (0,1).
\end{equation}

Now, since $\sup_{t\in [s,T]}|\Phi^v_{s,T}(t)|<\infty$, by Girsanov's theorem,
$$W_t^\vv:= W_t -\vv \int_{s}^t \si_r^*(\si_r\si_r^*)^{-1} \Phi^v_{s,T}(r) \d r,\ \ t\in [s,T]$$ is a cylindrical Brownian motion
on $\H_2$ under the probability measure $\d \P_\vv:= R_\vv\d\P$, where
\beq\label{2.7} R_\vv:= \exp\bigg[\vv \int_{s}^{T} \big\<\si_r^*(\si_r\si_r^*)^{-1} \Phi^v_{s,T}(r), \d W_r\big\>
-\ff{\vv^2} 2  \int_{s}^{T} \big|\si_r^*(\si_r\si_r^*)^{-1}\Phi^v_{s,T}(r)\big|^2\d r\bigg].\end{equation}
Hence, if we write \eqref{2.3} as
$$\beg{cases}
\d X_{s,t}^\vv= \big\{A_1X_{s,t}^\vv +BY_{s,t}^\vv\big\}\d t, & X^\vv_{s,s}=x+\vv v_1,\\
\d Y_{s,t}^\vv=  A_2Y_{s,t}^\vv\,\d t +\si_t\d W_t^\vv, &    Y_{s,s}^\vv= y+\vv v_2,
\end{cases}$$
then by the weak uniqueness of the solution, we obtain
$$
P_{s,T}^0 f(x+\vv v_1, y+\vv v_2)= \E_{\P_\vv} \big[f(X_{s,T}^\vv, Y_{s,T}^\vv)\big] =\E\big[ R_\vv f(X_{s,T}^\vv, Y_{s,T}^\vv)\big].
$$
Combining this with \eqref{2.6} and \eqref{2.7}, we arrive at
\beg{equation*}\beg{split} (\nn_v P_{s,T}^0 f)(x,y)& =\lim_{\vv\downarrow 0} \E\bigg[\ff{R_\vv-1}\vv f\big((X_{s,T}^0, Y_{s,T}^0)(x,y)\big)\bigg]\\
&=   \E\bigg[ f\big((X_{s,T}^0, Y_{s,T}^0)(x,y)\big)\int_{s}^{T} \big\<\si_r^*(\si_r\si_r^*)^{-1} \Phi^v_{s,T}(r), \d W_r\big\>\bigg].
\end{split}\end{equation*}
The proof is finished.
\end{proof}

\paragraph{Remark 2.1.}  In formula \eqref{BS}, although the operator $A_1$ does not appear explicitly, it is used in \eqref{2.5} implicitly though assumption {\bf (H2)}.
Of course, in the finite-dimensional case this assumption is not needed, see \cite{GW,Zh}.

\

We note that the derivative formula in Theorem \ref{T2.1} also applies to Hilbert-valued map $f\in \B_b(\H;\tt\H)$ by expanding
$f$ along an orthonormal basis of $\tt\H$, where $\tt\H$ is a  separable Hilbert space.
Moreover, by the semigroup property, formula \eqref{BS} also implies high order derivative formulas. For instance, for $t\in (s,T)$ and $v,\tt v\in \H$,   \eqref{2.2} implies
\begin{align*}
v_t:= \nn_{  v}(X_{s,t}^0, Y_{s,t}^0)&=\bigg(\e^{(t-s)A_1} v_1 +\int_{s}^t\e^{(t-r)A_1}B\e^{(r-s)A_2}v_2\,\d r,\ \e^{(t-s)A_2}v_2\bigg).
\end{align*}
Then by $P_{s,T}^0 f=P_{s,t}^0 P_{t,T}^0 f$ and  \eqref{2.2}, \eqref{BS}, we have
\beq\label{BS2} \beg{split} \nn_{v}\nn_{\tt v} P_{s,T}^0 f  &= \nn_v
\E\bigg[ (P_{t,T}^0f)(X_{s,t}^0, Y_{s,t}^0)\int_{s}^{t} \big\<\si_r^*(\si_r\si_r^*)^{-1} \Phi^{\tt v}_{s,t}(r), \d W_r\big\>\bigg]\\
 &=\E\bigg[ \big(\nn_{v}P_{t,T}^0f(X_{s,t}^0, Y_{s,t}^0)\big)\int_{s}^{t} \big\<\si_r^*(\si_r\si_r^*)^{-1} \Phi^{\tt v}_{s,t}(r), \d W_r\big\>\bigg]\\
 &= \E\bigg[ \big(\nn_{v_t}P_{t,T}^0f\big)(X_{s,t}^0, Y_{s,t}^0)
 \int_{s}^{t} \big\<\si_r^*(\si_r\si_r^*)^{-1} \Phi^{\tt v}_{s,t}(r), \d W_r\big\>\bigg]\\
 &= \E\bigg[f(X_{s,T}^0, Y_{s,T}^0)\bigg(\int_t^{T}\big\<\si_r^*(\si_r\si_r^*)^{-1} \Phi^{v_t}_{t,T}(r), \d W_r\big\>\bigg) \\
  &\qquad \qquad\qquad\qquad\times   \bigg(\int_{s}^{t} \big\<\si_r^*(\si_r\si_r^*)^{-1} \Phi^{\tt v}_{s,t}(r), \d W_r\big\>\bigg)\bigg].\end{split}\end{equation}

We will use \eqref{BS} and \eqref{BS2} to estimate derivatives of $P_{s,T}^0f$ for $f\in \B_b(\H;\tt\H)$ in terms of the norm
$$
\|f\|_{\phi,\psi}:=\|f\|_\infty+ \sup_{(x,y)\ne(x',y')\in\H} \ff{|f(x,y)-f(x',y')|}{\phi(|x-x'|)  + \psi(|y-y'|)},$$ where $\phi,\psi\in \D_0 $
and $\|\cdot\|_\infty$ is the uniform norm. Let
$$
\C_{\phi,\psi}(\H;\tt\H):=\Big\{f\in \B_b(\H;\tt\H):\ \|f\|_{\phi,\psi}<\infty\Big\}.
$$
Then $(\C_{\phi,\psi}(\H;\tt\H),\|\cdot\|_{{\phi,\psi}})$ is a  Banach space.
In particular, for any $\alpha\in [0,1]$, if we let
$\gamma_\aa(s)= s^\aa 1_{(0,\infty)}(s)$, then for $\alpha,\beta\in [0,1]$, $\C_{\gamma_\alpha,\gamma_\beta}(\H;\tt\H)$ is the usual H\"older space and
$$
\|f\|_{\gamma_\aa,\gamma_\beta}=\|f\|_\infty+ \sup_{(x,y)\ne(x',y')\in\H} \ff{|f(x,y)-f(x',y')|}{ |x-x'|^\aa   +  |y-y'|^\beta}.
$$
Note that $\|f\|_{\gamma_0,\gamma_0}\approx\|f\|_\infty$.
\beg{cor}\label{C2.2} Assume {\bf (H1)}-{\bf (H3)} and let $T>0$ be fixed. Let $\nn^{(i)}$ denote the gradient operator on $\H_i, i=1,2.$
\beg{enumerate}\item[$(1)$] There exists a constant $C>0$ such that for any $\aa\in [0,1],$
$$\|\nn^{(1)} P_{s,t}^0 f\|_\infty \le \ff{C\|f\|_{\gg_\aa,\gg_0}}{(t-s)^{\ff{3(1-\alpha)}2}},\  \ 0\le s<t\le T,
f\in \C_{\gamma_\alpha,\gamma_0}(\H;\tt\H).
$$
\item[$(2)$] There exists a constant $C>0$ such that for any $\aa\in [0,1]$ and $\phi\in\D_0$ with $\phi^2$  concave,
$$\|\nn^{(2)} P_{s,t}^0 f\|_\infty \le \ff{C\|f\|_{\gg_\aa,\phi}}{\ss{t-s}}\Big[(t-s)^{\ff{\aa (2+\delta)}2} +\phi\big(C(t-s)^{\ff{\delta}2}\big)\Big]$$
holds for all $0\le s<t\le T$ and
$ f\in \C_{\gamma_\alpha,\phi}(\H;\tt\H),$
where $\delta\in (0,1)$ is in {\bf (H3)}. In particular,
$$\|\nn^{(2)} P_{s,t}^0 f\|_\infty \le \ff{C\|f\|_{\infty}}{\ss{t-s}},\ \ 0\le s<t\le T, f\in \B_b(\H;\tt\H).$$\end{enumerate} \end{cor}
\beg{proof} (1) By the interpolation theorem (cf. \cite[Theorem 1.2.1]{Lu}), it suffices to prove it for $\aa=0,1.$

(1a) Let $\aa=1$. For any $v_1\in \H_1$, \eqref{2.2} implies
\beq\label{2.8}
\nn^{(1)}_{v_1} Y_{s,t}^0=0,\ \ \nn^{(1)}_{v_1} X_{s,t}^0 =\e^{(t-s)A_1}v_1.
\end{equation}
So, for any $f\in C_b^1(\H;\tt\H)$,
\beg{equation*}\beg{split}
|\nn_{v_1}^{(1)} P_{s,t}^0f|&=\left|\E\left[\left(\nn^{(1)}_{\nn_{v_1}^{(1)}X_{s,t}^0}f\right)(X_{s,t}^0, Y_{s,t}^0)\right]\right|
\le \|f\|_{1,0}|\e^{(t-s)A_1}v_1|\le \|f\|_{1,0}|v_1|.
\end{split}\end{equation*}
Thus, assertion (1) is proved for $\aa=1.$

(1b) Let $\aa=0$ and $0\le s<t\le T, v\in\H$. By {\bf (H1)} and the definitions of $\Phi^v_{s,t}$ and $V^v_{s,t}$, there exists a constant $C_1>0$ such that
\beq\label{2.9} \beg{split}
\int_{s}^{t} |\si_r^*(\si_r\si_r^*)^{-1}\Phi^v_{s,t}(r)|^2\d r
  \le C_1\left[\ff{|v_1|^2}{(t-s)^3} + \ff{|v_2|^2}{t-s}\right].
  \end{split} \end{equation}
Combining this with \eqref{BS} we obtain
 \beq\label{2.10} \beg{split}
 |\nn_vP_{s,t}^0f|^2 &\le (P_{s,t}^0|f|^2)\int_{s}^{t} |\si_r^*(\si_r\si_r^*)^{-1} \Phi^v_{s,t}(r)|^2\d r\\
 &\le C_1(P_{s,t}^0|f|^2)\left[\ff{|v_1|^2}{(t-s)^3} + \ff{|v_2|^2}{t-s}\right].
 \end{split} \end{equation}
 In particular, with $v_2=0$ this implies assertion (1) for $\aa=0.$

 (2) For $v_2\in\H_2$ and $(x,y)\in\H$, let $(X_{s,t}^0,Y_{s,t}^0)= (X_{s,t}^0,Y_{s,t}^0)(x,y)$ and
$$\tt x= \e^{(t-s)A_1}x+ \int_{s}^{t} \e^{(t-r)A_1}B\e^{(r-s)A_2}y\d r,\ \ \tt y= \e^{(t-s)A_2}y.$$ Moreover, let
$$\xi= \int_{s}^{t} \e^{(t-r)A_1}B\d r\int_{s}^r \e^{(r-r')A_2} \si_{r'}\d W_{r'},\ \ \eta=\int_{s}^{t} \e^{(t-r)A_2} \si_r\d W_r.$$
Since $\E \int_{s}^{t} \big\<\si_r^*(\si_r\si_r^*)^{-1}\Phi^v_{s,t}(r), \d W_r\big\>=0$, applying \eqref{BS} with $v=(0,v_2)$ and using \eqref{2.9}, we obtain
\beq\label{2.11}\beg{split} |\nn_{v_2}^{(2)} P_{s,t}^0f|(x,y)
&\le \bigg|\E\bigg[\Big\{f(X_{s,t}^0, Y_{s,t}^0)-f(\tt x,\tt y)\Big\}\int_{s}^{t} \big\<\si_r^*(\si_r\si_r^*)^{-1} \Phi^v_{s,t}(r), \d W_r\big\>\bigg]\bigg|\\
&\le \|f\|_{\gg_\aa,\phi} \E\bigg[\Big(|\xi|^\aa +\phi(|\eta|)\Big)\bigg|\int_{s}^{t} \big\<\si_r^*(\si_r\si_r^*)^{-1} \Phi^v_{s,t}(r), \d W_r\big\>\bigg|\bigg]\\
&\le \ff{C\|f\|_{\gg_\aa,\phi}|v_2|}{\ss{t-s}} \ss{\E(|\xi|^{2\aa}+\phi(|\eta|))^2}.\end{split}\end{equation}
Noting that {\bf (H3)} implies
\beg{equation}\label{EW1}\beg{split}
& \int_{s}^{t} \|\e^{(t-r)A_2}\si_r\|_{HS}^2\d r \le c_1\sum_{i\ge 1} \int_{s}^{t} \e^{-2\ll_i(t-r)}\d r\\
&\le c_1 \sum_{i\ge 1} \ff{1-\e^{-2\ll_i(t-s)}}{2\ll_i} \le c_1\sum_{i\ge 1} \ff{(2\ll_i(t-s))^{\delta}}{2\ll_i} =c_2 (t-s)^{\delta},
\end{split}\end{equation}
for $c_1:= \sup_{t\in [0,T]}\|\si_t\|^2,
c_2:= 2^{\delta-1}c_1\sum_{i\ge 1} \ff 1{\ll_i^{1-\delta}}<\infty,$
by Jensen's inequality we have
$$
\E|\xi|^{2\aa}\le c_3(t-s)^{(2+\delta)\aa},\ \ \E\phi(|\eta|)^2\le (\phi(\E|\eta|))^2\le \big(\phi\big( [c_2(t-s)]^{\delta/2}\big)\big)^2$$
for some constant $c_3>0.$ Combining this with \eqref{2.11} we prove the first assertion in (2), which implies the second assertion by taking $\aa=0$ and $\phi=\gamma_0=1.$
\end{proof}

\beg{cor}\label{C2.3} Assume {\bf (H1)}-{\bf (H3)}. For any $T>0, \aa\in [0,1]$ and $\phi\in\D_0$ with concave $\phi^2$, there exist constants
$C_1, C_2, C_3>0$ such that for any $f\in \B_b(\H;\tt\H),$
\beg{equation*}\beg{split}
&\|\nn_{\tt v_2}^{(2)}\nn_{v_2}^{(2)} P_{s,t}^0f\|_\infty\le
C_1\|f\|_{\gg_\aa,\phi}|v_2| |\tt v_2|
 \ff{(t-s)^{\aa(2+\delta)/2}+\phi(C_2(t-s)^{\delta/2})}{t-s},\\
&\|\nn_{v_1}^{(1)}\nn_{v_2}^{(2)} P_{s,t}^0f\|_\infty\le
\ff{C_3\|f\|_{\gg_\aa,\gg_0}|v_1| |v_2|}{(t-s)^{(4-3\aa)/2}},\ \ 0\le s<t\le T, v_1\in \H_1, v_2,\tt v_2\in \H_2.
\end{split}\end{equation*}
\end{cor}

\beg{proof} By the second equality in \eqref{BS2} for $v= (0,v_2)$ and $\tt v= (0,\tt v_2),$ we obtain
$$
\nn_{\tt v_2}^{(2)}\nn_{v_2}^{(2)} P_{s,t}^0f =\E\left[  \Big(\nn_{\tt v_2}^{(2)}P_{\frac{s+t}{2},t}^0f(X_{s,\frac{s+t}{2}}^0, Y_{s,\frac{s+t}{2}}^0)\Big)
\int_{s}^{\frac{s+t}{2}} \big\<\si_r^*(\si_r\si_r^*)^{-1}
\Phi^v_{s,\frac{s+t}{2}}(r), \d W_r\big\>\right].
$$
Combining this with the first inequality in Corollary \ref{C2.2} (2) and \eqref{2.9}, we derive
\beg{equation*}
\beg{split}  \|\nn_{\tt v_2}^{(2)}\nn_{v_2}^{(2)} P_{s,t}^0f\|_\infty
&\le\|\nn^{(2)}_{\tt v_2} P_{\frac{s+t}{2},t}^0f\|_\infty
\left[\E \int_{s}^{\frac{s+t}{2}} \big|\si_r^*(\si_r\si_r^*)^{-1} \Phi^v_{s,\frac{s+t}{2},v}(r)\big|^2\d r\right]^{\ff 1 2} \\
&\le C_1\|f\|_{\gg_\aa,\phi} |\tt v_2||v_2|
\ff{(t-s)^{\aa(2+\delta)/2}+\phi(C_2(t-s)^{\delta/2})}{t-s}.
\end{split}
\end{equation*}

Similarly, with $\tt v= (v_1,0)$ in place of $(0,\tt v_2)$ the second equality in \eqref{BS2} implies
$$
\nn_{v_1}^{(1)}\nn_{v_2}^{(2)} P_{s,t}^0f=\E\left[  \Big(\nn_{v_1}^{(1)}P_{\frac{s+t}{2},t}^0f\big(X_{s,\frac{s+t}{2}}^0, Y_{s,\frac{s+t}{2}}^0\big)\Big)
\int_{s}^{\frac{s+t}{2}} \big\<\si_r^*(\si_r\si_r^*)^{-1} \Phi^v_{s,\frac{s+t}{2}}(r), \d W_r\big\>\right].
$$
By using Corollary \ref{C2.2} (1) and \eqref{2.9}, we prove the second inequality.
\end{proof}

Finally,  we apply the above derivative estimates to the resolvent
$$R_{s,t}^\ll f:= \int_{s}^{t} \e^{-(r-s)\ll} P_{s,r}^0 f_r\d r,\ \ \ll\ge 0, 0\le s\le t, f\in \B_b([0,T]\times\H;\tt\H),$$
which will be used in the next section to construct the regularization transform. For any $f\in \B_b([0,T]\times\H;\tt\H)$, we simply denote
$$
\|f\|_{\phi,\psi}= \sup_{t\in [0,T]} \|f_t\|_{\phi,\psi},\ \ \ \phi,\psi\in\D_0.
$$

\beg{cor}\label{C2.4} Assume {\bf (H1)}-{\bf (H3)} and let $T>0$ be fixed.
 \beg{enumerate} \item[$(1)$] $R_{s,t}^\ll f\in \B_b([0,T]; \scr C_{\gamma_0,\gamma_1})$ for any $\ll\ge 0, 0\le s\le t$ and $f\in \B_b([0,T]\times\H;\tt\H).$
 Moreover, there exists a constant $C>0$ such that for any $f\in \B_b([0,T]\times\H;\tt\H),$
 $$\|\nn^{(2)} R_\cdot^\ll f\|_\infty:= \sup_{0\le s\le t\le T, z\in\H} \|\nn^{(2)} R_{s,t}^\ll f(z)\|\le C\|f\|_{\infty}\left[\ff{1-\e^{-\ll T}}\ll\right]^{\ff 1 2},\ \ \ll\ge 0.$$
 Finally, for any $\aa\in (\ff 2 3,1)$ there exists a  function $\theta: [0,\infty)\to (0,\infty)$ with $\theta(\ll)\downarrow 0$ as $\ll\uparrow\infty$ such that
 $$\|\nn^{(1)} R_\cdot^\ll\|_\infty \le \theta(\ll) \|f\|_{\gg_\aa,\gg_0},\ \ \ll\ge 0, f\in \B_b([0,T]; \scr C_{\gg_\aa,\gg_0}(\H;\tt\H)).$$
\item[$(2)$] For any $\aa\in (\ff 2 3,1)$ and $\phi\in \D_1$,   there exists $0<\theta(\ll)\downarrow 0$ as $\ll\uparrow\infty$ such that
 $$\|\nn\nn^{(2)}R^\ll_{s,t}f\|_\infty \le \theta(\ll) \|f\|_{\gg_\aa,\phi},\ \ 0\le s\le t\le T, f\in
 \B_b([0,T]; \scr C_{\gg_\aa,\phi}(\H;\tt\H)).
 $$

\item[$(3)$] Let $\aa\in (\ff 2 3,1)$ and $\phi\in\D_0$ with concave $\phi^2$.
Then there exists a constant $C>0$ such that for any $f\in\B_b([0,T]\times\H;\tt\H),$
 $$
 \|\nn^{(2)}(R_{s,t}^\ll f)(z)- \nn^{(2)}(R_{s,t}^\ll f)(z')\|\le C\|f\|_{\gamma_\aa,\phi}
 \inf_{r\in (0,1)} \bigg\{r+|z-z'|\bigg(1+\int_{r^{\delta}}^1 \ff{\phi(s)} s\,\d s\bigg)\bigg\}
 $$ holds for all $\ll\ge 0$ and $0\le s\le t\le T,$ where $\delta\in (0,1)$ is in {\bf (H3)}.
 \end{enumerate}\end{cor}

 \beg{proof} Assertion (1) follows immediately from the definition of $R^\ll_{s,t}$ and Corollary \ref{C2.2}.   Next,
since $\|\cdot\|_{\gg_\aa,\gg_0}\le c\|\cdot\|_{\gg_\aa,\phi}$ holds for some constant $c>0$, Corollary \ref{C2.3} implies
 \beq\label{W0}
 \|\nabla\nn^{(2)} P_{s,r}^0 f\|_\infty\le C\left[\ff 1 {(r-s)^{(4-3\alpha)/2}}
 + \ff{\phi(C(r-s)^{\delta/2})}{r-s}\right]\|f\|_{\gamma_\aa,\phi}.\end{equation}
Hence,
\beq\label{NW}\beg{split} \|\nabla\nn^{(2)}(R_{s,t}^\ll f)\|_\infty
\le C \|f\|_{\gamma_\aa,\phi}\int_{s}^{t} \e^{-\ll (r-s)}\left[\ff{1}{(r-s)^{(4-3\alpha)/2}} +\ff{\phi(C(r-s)^{\delta/2})}{r-s}\right]\d r,
 \end{split} \end{equation}
 which implies the assertion (2) due to $\aa\in (\ff 2 3,1)$ and $\int^1_0\phi(s)/s\d s<\infty$.

 Finally, we prove (3) without assuming $\int^1_0\phi(s)/s\d s<\infty$. By Corollary 2.2 (2), we have
 $$
 \|\nn^{(2)} P_{s,r}^0f(z)- \nn^{(2)} P_{s,r} f(z')\|\le
 \ff{C\|f\|_{\gamma_\aa,\phi}}{\ss{r-s}},\ \ r>s.
 $$
 On the other hand, by \eqref{W0}, we have
 \beq\label{W00}
 |\nn^{(2)} P_{s,r}^0 f(z)-\nn^{(2)} P_{s,r}^0 f(z')|\le C|z-z'|\left[
 \ff 1 {(r-s)^{\ff{4-3\aa} 2}}+ \ff{\phi(C(r-s)^{\ff{\delta}2})}{r-s}\right]\|f\|_{\gamma_\aa,\phi}.
 \end{equation}
 Combining these together we obtain
\beq\label{NW1}\beg{split}
&\|\nn^{(2)}(R_{s,t}^\ll f)(z)- \nn^{(2)}(R_{s,t}^\ll f)(z')\|\\
 &\le C \|f\|_{\gamma_\aa,\phi}\int_{s}^{t} \e^{-\ll (r-s)}\left[\ff{|z-z'|}{(r-s)^{\ff{4-3\aa}2}} +\ff 1 {\ss{r-s}}\land \ff{|z-z'|\phi(C(r-s)^{\ff{\delta}2})}{r-s}\right]\d s\end{split} \end{equation} for some constant $C>0.$ Noting that
 \beg{equation*}\beg{split} &\int_0^T \e^{-\ll s} \left[\ff 1 {\ss s}\land \ff{|z-z'|\phi(Cs^{\ff{\delta} 2})} s\right]\d s \\
 &\le \int_0^{r^2 /C^{2/\delta}} \ff 1 {\ss s}\d s +|z-z'| \int_{(r^2 /C^{2/\delta})\land T}^T \ff{\phi(Cs^{\ff{\delta}2})} s \,\d s\\
 &\le C' r+  C' |z-z'|\bigg(1+ \int_{r^{\delta}}^1 \ff{\phi(s)} s \,\d s\bigg),\ \ r\in (0,1)\end{split}\end{equation*} holds for some constant $C'>0$, the assertion (3) follows from \eqref{NW}.
 \end{proof}

\section{Regularization transform}

Throughout this section, we assume that $b$ is bounded. We aim to construct a map $\Theta: [0,T]\times \H\to\H$ such that
$\Theta_t$ is a diffeomorphism for every $t\in [0,T]$, and $(\bar X_t, \bar Y_t):=\Theta_t(X_t,Y_t)$ solves
an equation which has pathwise uniqueness
provided $(X_t,Y_t)$ solves \eqref{1.1}. In this way we prove the pathwise uniqueness of \eqref{1.1}.

The transform will be constructed as
\beq\label{3.0} \Theta_s(x,y):= (x, y+ u_s^\ll(x,y)),\ \ s\in [0,T], (x,y)\in\H\end{equation}  for large $\ll>0$,
where $u_s^\ll$ solves
the integral equation
\beq\label{3.1} u_s^\ll= \int_s^T \e^{-\ll(t-s)}P_{s,t}^0 \Big\{\nn_{b_t}^{(2)} u_t^\ll+b_t\Big\}\d t,\ \ s\in [0,T].\end{equation}
To ensure that $\Theta_s$ is a bijection on $\H$ for every $s\in [0,T]$, we need $\|\nn^{(2)}u_s^\ll\|_\infty<1.$ So, we first solve \eqref{3.1} and estimate  $\nn^{(2)}u_s^\ll.$

\subsection{Gradient estimates on $u_s^\ll$}

\beg{prp}\label{P3.1} Assume  {\bf (H1)}-{\bf (H3)} and let $\aa\in (\ff 2 3,1), \phi\in \D_0$ and $T>0$.
\beg{enumerate}
\item[$(1)$] For any $R>0$ there exists a constant $\ll(R)>0$ such that $\eqref{3.1}$ has a unique solution
$u^\ll\in C_b([0,T];\C_{\gamma_0,\gamma_1}(\H;\H_2))$ provided $\|b\|_\infty\le R$.
\item[$(2)$] For any $R>0$ there exist a constant $\ll(R)>0$ and a positive function $\theta$ on $[\ll(R),\infty)$ with $\theta(\ll)\downarrow 0$ as $\ll\uparrow\infty$
such that if $\|b\|_{\gamma_\aa,\phi}\le R$ then
\beq\label{3.2'} \sup_{s\in [0,T]} \|\nn u^\ll_s\|_\infty\le \dd(\ll)\end{equation} and
\beq\label{3.2}  \sup_{s\in [0,T]}\|\nn^{(2)}u_s^\ll(z)- \nn^{(2)}u_s^\ll(z')\|\le C\inf_{r\in (0,1)} \bigg\{r+ |z-z'|\bigg(1+ \int_{r^{\delta}}^1 \ff{\phi(s)}s\,\d s\bigg)\bigg\}
 \end{equation} hold  for all $ \ll\ge \ll(R)$ and $z,z'\in \H$,  where $\delta\in (0,1)$ is in {\bf (H3)}.
\item[$(3)$] If $\phi\in\D_1$, then for any $R>0$,
$$\lim_{\ll\to\infty} \sup_{t\in [0,T]} \|\nn \nn^{(2)}u_t^\ll\|_\infty=0$$ holds uniformly for $b$ with $\|b\|_{\gamma_\aa,\phi}\le R$. \end{enumerate}
\end{prp}

\beg{proof} (1) Let
$$\scr H= C_b([0,T];\C_{\gamma_0,\gamma_1}(\H;\H_2)),$$
which is a Banach space with the norm
$$\|f\|_{\scr H} := \sup_{t\in [0,T]} \|f_t\|_{\C_{\gamma_0,\gamma_1}}=\sup_{t\in [0,T]}\big(\|f_t\|_{\gamma_0,\gamma_1}+\|f_t\|_\infty\big).$$
For any $f\in\scr H$, let
$$
\GG_s^\ll(f):= R_{s,T}^\ll (\nn_b^{(2)}f+b)= \int_s^T \e^{-\ll(t-s)}P_{s,t}^0 \{\nn_{b_t}^{(2)}f_t +b_t\}\d t,\ \ s\in [0,T].
$$
By Corollary \ref{C2.4} (1), there exists a constant $\ll(R)>0$ such that  $\|b\|_\infty\le R$ implies
\beq\label{W1} \|\GG^\ll(f)-\GG^\ll(g)\|_{\scr H}\le \ff 1 2 \|f-g\|_{\scr H},\ \ f,g\in \scr H,\ll\ge \ll(R).\end{equation}
Thus, by the fixed-point theorem, for $\ll\ge \ll(R)$ and $\|b\|_\infty\le R$ the equation \eqref{3.1} has a unique solution
$u^\ll\in {\scr H}.$    Moreover, taking $g=0$ in \eqref{W1} and by Corollary \ref{C2.4} (1), we obtain
\beq\label{W2}
\|u^\ll\|_{\scr H}\le 2\|\GG^\ll(0)\|_{\scr H} \le \ff{C\|b\|_\infty}{\sqrt{\ll}},\ \ \ll\ge \ll(R).
\end{equation}
In particular,
\beq\label{WL1}
\|\nn^{(2)} u^\ll\|_\infty\le \ff{C_1}{\ss\ll},\ \ \ll\ge \ll(R).
\end{equation}

(2)   All constants $C_i$ mentioned in the proof below are uniform in $b$ with $\|b\|_{\gamma_\aa,\phi}\le R$ and $\ll\ge\ll(R)$.
By \eqref{W0} and \eqref{W2}, we have
\beg{equation*}\beg{split} &\|\nn^{(2)}P_{s,t}^0\{\nn_{b_t}^{(2)}u_t^\ll+b_t\}(z) - \nn^{(2)}P_{s,t}^0\{\nn^{(2)}_{b_t}u_t^\ll+b_t\}(z')\|\\
&\qquad\le C_2\|u^\ll\|_{\scr H} \ff{|z-z'|}{(t-s)^2} \le \ff{C_3|z-z'|}{(t-s)^2},\ \ 0\le s<t\le T.\end{split}\end{equation*}
Combining these with \eqref{3.1} we obtain
$$ \|\nn^{(2)}u^\ll(z)- \nn^{(2)}u^\ll(z')\|\le (C_1\lor C_3)\int_0^T \Big(1\land \ff{|z-z'|}{t^2}\Big)\d t\le C_4 \ss{|z-z'|}.$$
Since $\|b\|_{\gamma_\aa,\gamma_0}\le \|b\|_{\gamma_\aa,\phi}\le R$, this implies
\beq\label{Y*}\|\nn^{(2)}_{b}u^\ll+b\|_{\gamma_{1/2},\gamma_0}\le C_5.\end{equation}
Combining this with \eqref{W0} for $\aa=\ff 1 2$, we obtain
$$\|\nn^{(2)}P_{s,t}^0\{\nn_{b_t}u_t^\ll+b_t\}(z) - \nn^{(2)}P_{s,t}^0\{\nn_{b_t}u_t^\ll+b_t\}(z')\|\le \ff{C_6|z-z'|}{(t-s)^{5/4}},\ \ 0\le s<t\le T.$$
This together with \eqref{WL1} leads
$$ \|\nn^{(2)}u^\ll(z)- \nn^{(2)}u^\ll(z')\|\le (C_1\lor C_6)\int_0^T \Big(1\land \ff{|z-z'|}{t^{5/4}}\Big)\d t\le C_7  |z-z'|^{4/5}.$$
Let $\psi(s)=\ss{\phi(s)^2+s}$ such that $\psi\in\D_0$.
Since $\|b\|_{\gamma_\aa,\psi}\le \|b\|_{\gamma_\aa,\phi}\le R$, this implies
$$\|\nn_b^{(2)}u^\ll+b\|_{\gamma_{\aa\land \ff 4 5},\psi} \le C_8.$$
Thus, by Corollary \ref{C2.4} (3) for $(\aa\land \ff 4 5,\psi)$ in place of $(\aa,\phi)$, we obtain
\beg{equation*}\beg{split}
\|\nn^{(2)}u_s^\ll(z)- \nn^{(2)}u_s^\ll(z')\|&\le C_9\inf_{r\in (0,1)} \bigg\{r+ |z-z'|\bigg(1+ \int_{r^{\delta}}^1 \ff{\psi(s)}s\,\d s\bigg)\bigg\}\\
&\le C_{10} \inf_{r\in (0,1)} \bigg\{r+ |z-z'|\bigg(1+ \int_{r^{\delta}}^1 \ff{\phi(s)}s\,\d s\bigg)\bigg\}.
\end{split}\end{equation*}
Therefore, \eqref{3.2} holds for some constant $C>0.$
Finally, by \eqref{3.1}, \eqref{Y*} and the last inequality in Corollary \ref{C2.4} (1),
we prove \eqref{3.2'} for some positive function $\dd$ with $\dd(\ll)\downarrow 0$ as $\ll\uparrow \infty$. So, the proof of assertion (2) is finished.

(3) Finally, if $\|b\|_{\gamma_\alpha,\phi}\le R$ for some $\aa\in (\ff 2 3,1)$ and $\phi\in\D_1$,
then by assertion (2), $\nn^{(2)}u^\ll_s$ is Lipschitz continuous  uniformly in $s\in [0,T]$ so that
$$\|\nn_b^{(2)}u^\ll+b\|_{\gamma_\alpha,\phi} \le C_{11}.$$ So, by Corollary \ref{C2.4} (2) and \eqref{3.1}, we prove assertion (3).
\end{proof}

According to \eqref{W2} we have $\|\nn^{(2)}u^\ll\|_\infty<1$ for large $\ll>0$, so that $\Theta_s$ given in \eqref{3.0} is a bijection for every $s\in [0,T]$. To figure out the equation satisfied by
$(\bar X_t, \bar Y_t):=\Theta_t(X_t,Y_t),$ we need to apply It\^o's formula to $u_t^\ll(X_t,Y_t)$, which is however not available in the infinite-dimensional setting. Thus, we need to investigate approximations on $u^\ll$ such that It\^o's formula can be established for $u_t^\ll(X_t,Y_t)$ by finite-dimensional approximations.

\subsection{Approximations on $u^\ll$}

Let $\H_i^{(n)} (i=1,2), \H^{(n)}$ and $\pi^{(n)}=(\pi_1^{(n)},\pi_2^{(n)})$ be defined in Section 1. Since
$BB^*$
is invertible, for any $x\in \H_1$ we have
$x=By'$ for
$y'= B^*(BB^*)^{-1} x\in \H_2.$
So,
\beq\label{P1} \lim_{n\to\infty} |\pi^{(n)}(x,y)-(x,y)|\le \lim_{n\to\infty} \big(|\pi_2^{(n)}y-y|+ |B(\pi_2^{(n)} y'-y')|\big)=0,\ \ (x,y)\in\H.\end{equation}
Let
$$
b_s^{(n)} =\pi_2^{(n)} b_s, \ \ \si_s^{(n)}= \pi_2^{(n)} \si_s,\ \ s\ge 0.
$$
Since $\H_2^{(n)}$ is $A_2$-invariant, $\H_1^{(n)}=B\H_2^{(n)}$ and $A_1B=BA_2$, we have
\beq\label{P2}
\pi_2^{(n)} A_2\pi_2^{(n)}= A_2\pi_2^{(n)},\ \ \pi_1^{(n)} A_1\pi_1^{(n)}= A_1 \pi_1^{(n)}.
\end{equation}
For any $n\ge 1$, consider the equation
\beq\label{P4}
u_s^{\ll,n} = \int_s^T \e^{-\ll(t-s)}P_{s,t}^0\Big\{\nn_{b_t^{(n)}}^{(2)}u_t^{\ll,n} +b_t^{(n)}\Big\}\d t,\ \ s\in [0,T].
\end{equation}

\beg{prp}\label{P3.2} Assume  {\bf (H1)}-{\bf (H4)} and let $n_0$ be from {\bf (H4)} and $b\in \B_b([0,T]; C_b(\H;\H_2))$. Then
there exists a constant $\ll_0>0$ such that for any $\ll\ge \ll_0$ and $n\ge n_0$,
equation $\eqref{P4}$ has a unique solution $u^{\ll, n}\in C_b([0,T]; \C_{\gamma_0,\gamma_1}(\H;\H_2^{(n)}))$ with
\beq\label{P23}
\lim\limits_{\ll\to\infty}\sup_{n\ge n_0} \big(\|\nn^{(2)}u^{\ll,n}\|_\infty+ \|u^{\ll,n}\|_\infty\big)=0.
\end{equation}
Moreover, for all  $s\in[0,T]$ and  $z\in\H$, we have
\beq\label{P22}
\lim\limits_{n\to\infty} \big(|u^{\ll,n}_s(z)- u^\ll_s(z)|+ \|\nn^{(2)}(u^{\ll,n}_s-u^\ll_s)(z)\|\big)=0,\ \ll\ge \ll_0.
\end{equation}
\end{prp}

\beg{proof}  By Proposition \ref{P3.1} (1), there exists $\ll_0>0$ such that for any $n\ge 1$ and $\ll\ge \ll_0$, the equation \eqref{P4} has a unique solution $u^{\ll, n}\in C_b([0,T]; \C_{\gamma_0,\gamma_1}(\H;\H_2)).$ By {\bf (H4)} and \eqref{P2} we have
$$\pi^{(n)}_2 P_{s,t}^0 f= P_{s,t}^0 \pi_2^{(n)} f,\ \ n\ge n_0, f\in \B_b(\H;\H_2).$$ So, it is easy to see that $\pi_2^{(n)} u^{\ll,n}$ also solves the equation
\eqref{P4}. By the uniqueness we get $u^{\ll, n}\in C_b([0,T]; \C_{\gamma_0,\gamma_1}(\H;\H_2^{(n)}))$.
Obviously, \eqref{P23} follows from \eqref{W2}. Let us  prove \eqref{P22}.
By \eqref{3.1} and \eqref{P4}, we have
\beq\label{*W1}
\beg{split} u_s^\ll-u_s^{\ll,n} &= \int_s^T \e^{-\ll(t-s)}P_{s,t}^0\Big\{\nn_{b_t^{(n)}}^{(2)}(u_t^\ll-u_t^{\ll,n})\Big\}\d t\\
&\quad +\int_s^T \e^{-\ll(t-s)}P_{s,t}^0\Big\{\nn_{b_t-b_t^{(n)}}^{(2)}u_t^\ll +b_t-b_t^{(n)}\Big\}\d t.\end{split}
\end{equation}
Let
$$
g_s(z)= \limsup_{n\to\infty} \|\nn^{(2)} (u_s^\ll- u_s^{\ll,n})(z)\|,\ \ s\in [0,T], z\in\H.
$$
By \eqref{P23}, $g\in \B_b([0,T]\times\H).$ Moreover, since
$$
\|\nn^{(2)} u^\ll\|_\infty<\infty,\ \ \|b^{(n)}\|_\infty\le \|b\|_\infty<\infty,\ \ \lim_{n\to\infty} |b_t(z)-b_t^{(n)}(z)|=0,
$$
it follows from \eqref{2.10} with $v_1=0$ and the dominated convergence theorem that
\beg{equation*}\beg{split} &\limsup_{n\to\infty} \int_s^T \Big\|\nn^{(2)} P_{s,t}^0\Big\{\nn_{b_t-b_t^{(n)}}^{(2)} u_t^\ll + b_t -b_t^{(n)}\Big\}(z)\Big\|\d t\\
&\le C \limsup_{n\to\infty} \int_s^T \ff 1 {\ss{t-s}} \Big(P_{s,t}^0 \big|\nn_{b_t- b_t^{(n)}}^{(2)} u_t^\ll +b_t-b_t^{(n)}\big|^2 (z)\Big)^{\ff 1 2} \d t=0.\end{split}\end{equation*}
Combining this with \eqref{*W1} and \eqref{2.10}, we obtain
\beg{equation*}
\beg{split} g_s(z) &\le C \limsup_{n\to\infty} \int_s^T \ff {\e^{-\ll(t-s)}} {\ss{t-s}} \Big(P_{s,t}^0 \big|\nn_{ b_t^{(n)}}^{(2)} (u_t^\ll-u_t^{\ll,n})\big|^2(z)\Big)^{\ff 1 2}\d t \\
&\le C\|b\|_\infty \int_s^T \ff {\e^{-\ll(t-s)}} {\ss{t-s}} (P_{s,t}^0 |g|^2)^{\ff 1 2}(z)\d t\le \ff{C'}{\ss \ll} \|g\|_\infty,\ \ s\in [0,T], z\in\H, \ll\ge \ll_0\end{split}
\end{equation*}
for some constants $C,C'>0$ independent of $s,z$.
This implies $g_s(z)=0$ when $\ll$ is large enough.
\end{proof}

\section{Representation of $Y_t$ using $u^\ll$}
Assume {\bf (H1)}-{\bf (H4)} and let $b\in \B_b([0,T]; C_b(\H;\H_2))$ for some $T>0.$ Let $u^\ll$ be constructed in Section 3 for $\ll\ge \ll_0$, where $\ll_0 $ is in Proposition \ref{P3.2}.

\beg{thm}\label{T4.1} Assume {\bf (H1)}-{\bf (H4)} and let $b\in \B_b([0,T]; C_b(\H;\H_2))$ for some $T>0.$
If $(Z_t)_{t\in[0, T\land\tau]}:=(X_t,Y_t)_{t\in[0, T\land\tau]}$
solves $\eqref{1.1}$ for some stopping time $\tau$, i.e. $\P$-a.s.
\beq\label{4.1} \left\{\beg{aligned}&X_t= \e^{tA_1}X_0 +\int_0^t\e^{(t-s)A_1}BY_s\d s,\ 0\le t\le T\land\tau,\\
&Y_t =\e^{tA_2}Y_0 +\int_0^t \e^{(t-s)A_2}\big\{b_s(X_s,Y_s)\d s +\si_s \d W_s\big\},\end{aligned}\right.\end{equation}
then for any $\ll\ge \ll_0$, $\P$-a.s.
\beq\label{4.2} \beg{split}Y_t =\ &\e^{tA_2}Y_0 +\e^{tA_2}u_0^\ll(Z_0)- u_t^\ll(Z_t) +\int_0^t (\ll-A_2)\e^{(t-s)A_2}u_s^\ll(Z_s)\d s\\
&+ \int_0^t \e^{(t-s)A_2} \Big\{\si_s\d W_s+ \big(\nn^{(2)}_{\si_s\d W_s}u_s^\ll\big)(Z_s)\Big\},\ \ \ 0\le t\le T\land\tau.\end{split}\end{equation} \end{thm}

\beg{proof}
Let $\big(X_{s,t}^0, Y_{s,t}^0\big)_{t\ge s}$ solve the linear equation \eqref{2.1} with initial value $z\in\H$. Write
$$
Z^{0}_{s,t}:=\big(X_{s,t}^0, Y_{s,t}^0\big),\ \ Z^{(n)}_{s,t}:=\pi^{(n)}Z^0_{s,t}.
$$
By {\bf (H4)} and \eqref{P2}, $Z^{(n)}_{s,t}=\big(X_{s,t}^{(n)}, Y_{s,t}^{(n)}\big)$ solves the following equation:
\beq\label{4.3}
\beg{cases}
\d X_{s,t}^{(n)}= \big\{A_1X_{s,t}^{(n)} +BY_{s,t}^{(n)}\big\}\d t,\\
\d Y_{s,t}^{(n)}=  A_2Y_{s,t}^{(n)} \d t +\si_t^{(n)} \d W_t.
\end{cases}
\end{equation}
Since $Z_{s,s}^{(n)}(z)=\pi_n z$, by the uniqueness we have
\beq\label{43}
Z_{s,t}^{(n)}(z)= Z_{s,t}^{(n)}(\pi^{(n)} z).
\end{equation}
Let $u^\ll$ and $u^{\ll,n}$ be in Proposition \ref{P3.2}. Let
$$f_t^{(n)}= \nn^{(2)}_{b_t^{(n)}} u_t^{\ll,n} + b_t^{(n)},\ \ t\in [0,T], n\ge n_0.$$
For fixed $\vv_0>0$, let
$$
F_{s,r}^{(n)}(z):=P_{s,r+\vv_0}^0(f_r^{(n)}\circ\pi^{(n)})(z),\ \ 0\le s\le r\le T, z\in\H, n\ge n_0.
$$
By \eqref{43} we have
\beq\label{4.4}
F_{s,r}^{(n)}(z)= \E \Big[f_r^{(n)}\big(Z^{(n)}_{s,r+\vv_0}(z)\big)\Big]
= \E\Big[f_r^{(n)}\big(Z_{s,r+\vv_0}^{(n)}(\pi^{(n)} z)\big)\Big]=F_{s,r}^{(n)}(\pi^{(n)}(z)).
\end{equation}
Since $\vv_0>0$ is fixed, by \eqref{BS2} we have $F_{s,r}^{(n)}\in C_b^2(\H^{(n)};\H^{(n)}_2)$ with
\beq\label{4.5}
\sup_{0\le s\le r\le T, n\ge n_0} \Big(\|F_{s,r}^{(n)}\|_\infty + \|\nn F_{s,r}^{(n)}\|_\infty+ \|\nn\nn F_{s,r}^{(n)}\|_\infty\Big)<\infty.
\end{equation}
By \eqref{4.4}, \eqref{4.3} and It\^o's formula, we have
\beq\label{4.6}\beg{split}
\d F_{s,r}^{(n)}(Z_{s,t}^{(n)})
=(\scr L_t^{(n)} F_{s,r}^{(n)} )(Z_{s,t}^{(n)})\d t +\Big(\nn_{\si_t^{(n)}\d W_t}^{(2)} F_{s,r}^{(n)}\Big)(Z_{s,t}^{(n)}),
\end{split}\end{equation}
where for $F\in C_b^2(\H;\H_2)$ with $F=F\circ \pi^{(n)}$,
$$
(\scr L_t^{(n)}F)(x,y):=\ff 1 2 \sum_{i,j=1}^{n} \<\si_t\si_t^* e_i,e_j\>\nn_{e_i}^{(2)}\nn^{(2)}_{e_j} F(x,y)
+\nn_{A_1x+By}^{(1)} F(x,y),\ \ (x,y)\in\H.
$$
Combining \eqref{4.4}, \eqref{4.5} and \eqref{4.6}, we obtain
\beq\label{4.7}\beg{split}
\ff{\d}{\d s} F_{s,r}^{(n)} & =-\lim_{\vv\downarrow 0} \ff{F_{s-\vv,r}^{(n)}- F_{s,r}^{(n)}}\vv=-\lim_{\vv\downarrow 0}
\ff{\E [F_{s,r}^{(n)}(Z_{s-\vv,s}^{(n)})- F_{s,r}^{(n)}(Z_{s-\vv,s-\vv}^{(n)})]}\vv \\
&=- \lim_{\vv\downarrow 0} \ff 1 \vv \int_{s-\vv}^s \scr L_{s-\vv}^{(n)} F_{s,r}^{(n)}(Z_{s-\vv, t}^{(n)})\d t
 = - \scr L_s^{(n)} F_{s,r}^{(n)},\ \ \text{a.e.}\ s\in (0, \dd+r).
 \end{split}\end{equation}
Thus, if we write
$$
u_{s,\vv_0}^{\ll,n}:= \int_s^T \e^{-\ll(t-s)}P_{s,t+\vv_0}^0 (f_t^{(n)}\circ\pi^{(n)})\d t =\int_s^T \e^{-\ll(t-s)} F_{s,t}^{(n)}\d t,
$$
then $u_{s,\vv_0}^{\ll,n}$ satisfies
\beq\label{4.8}
\begin{split}
\pp_s u_{s,\vv_0}^{\ll,n}&=(\ll- \scr L_s^{(n)})u_{s,\vv_0}^{\ll,n}
-P_{s,s+\vv_0}^0 \big(\big\{\nn_{b_s^{(n)}}^{(2)} u_s^{\ll,n} + b_s^{(n)}\big\}\circ\pi^{(n)}\big).
\end{split}
\end{equation}
Noticing that by \eqref{4.4} and definitions,
$$
u_{s,\vv_0}^{\ll, n} = u_{s,\vv_0}^{\ll,n}\circ \pi^{(n)},
$$
we may apply It\^o's formula to $u_{s,\vv_0}^{\ll,n}(Z_s)= u_{s,\vv_0}^{\ll,n}(\pi^{(n)} Z_s)$ so that \eqref{4.1} and \eqref{4.8} imply
\beg{equation*}
\beg{split} &\d u_{s,\vv_0}^{\ll,n}(Z_s)  = \big(\nn_{\si_s\d W_s}^{(2)}\big)u_{s,\vv_0}^{\ll,n}(Z_s)
+\big(\nn_{b^{(n)}_s}^{(2)} u_{s,\vv_0}^{\ll,n}  +\scr L_s^{(n)} u_{s,\vv_0}^{\ll, n} +\pp_s u_{s,\vv_0}^{\ll, n}\big)(Z_s)\d s\\
&= \big(\nn_{\si_s\d W_s}^{(2)}\big)u_{s,\vv_0}^{\ll,n}(Z_s)+ \Big[\nn_{b_s^{(n)}}^{(2)} u_{s,\vv_0}^{\ll,n}-P_{s,s+\vv_0}^0
 \big(\big\{\nn_{b_s^{(n)}}^{(2)} u_s^{\ll,n} + b_s^{(n)}\big\}\circ\pi^{(n)}\big)
+\ll u_{s,\vv_0}^{\ll, n}\Big](Z_s)\d s.
\end{split}
\end{equation*}
Thus, for any $t\geq 0$, we have
\beq\label{4.9}\beg{split}
&u_{t,\vv_0}^{\ll, n}(Z_t)- \e^{tA_2} u_{0, \vv_0}^{\ll, n}(Z_0)=\int^t_0\d (\e^{(t-s)A_2} u_{s,\vv_0}^{\ll, n}(Z_s))/\d s\\
&=-\int_0^t  A_2  \e^{(t-s)A_2}   u_{s,\vv_0}^{\ll, n}(Z_s)\d s +\int_0^t \e^{(t-s)A_2} \big(\nn_{\si_s \d W_s}^{(2)} u_{s,\vv_0}^{\ll, n}\big)(Z_s)\\
&\quad+ \int_0^t \e^{(t-s)A_2} \Big(\nn_{b_s^{(n)}}^{(2)} u_{s,\vv_0}^{\ll, n}
- P_{s, s+\vv_0}^0 \big(\nn_{b_s^{(n)}}^{(2)} u_{s}^{\ll, n}\big)\Big)(\pi^{(n)}(Z_s))\d s\\
 &\quad+ \int_0^t \e^{(t-s)A_2}\Big(\ll u_{s,\vv_0}^{\ll, n}-P_{s,s+\vv_0}^0 b_s^{(n)}\Big)(\pi^{(n)}(Z_s))\d s.
 \end{split}
 \end{equation}
 We claim that the desired formula \eqref{4.2} follows by first letting $\vv_0\downarrow 0$ and then $n\uparrow \infty.$
To see this, we write
\beq\label{4.10} u_{s,\vv_0}^{\ll, n}- u_s^{\ll, n} =\int_s^T \e^{-\ll (t-s)}P_{s,t}^0 \big\{P_{t,t+\vv_0}^0 (f_t^{(n)}\circ\pi^{(n)} )- f_t^{(n)}\big\}\d t.\end{equation} By the boundedness and continuity of $b_s^{(n)}$ and $f_t^{(n)}$, and noting that $\lim_{\vv_0\downarrow 0} P_{t,t+\vv_0}^0 f=f$ for
$f\in C_b(\H;\H_2)$, we obtain
\beq\label{4.11}
\lim_{n\to\infty} \lim_{\vv_0\downarrow 0} (u_{s,\vv_0}^{\ll, n}- u_s^{\ll, n})=0,\ \ \lim_{n\to\infty} \lim_{\vv_0\downarrow 0}(P_{s, s+\vv_0}^0 b_s^{(n)})(\pi^{(n)}(Z_s))= b_s(Z_s).
\end{equation}
Moreover,   \eqref{2.10} with $v_1=0$ and \eqref{4.10} imply
$$\lim_{n\to\infty} \lim_{\vv_0\downarrow 0} \|\nn^{(2)}(u_{s,\vv_0}^{\ll, n}- u_s^{\ll, n})\|\le \lim_{n\to\infty} \lim_{\vv_0\downarrow 0}
\int_s^T \ff C{\ss{t-s}} \ss{P_{s,t}^0 |P_{t,t+\vv_0}^0 (f_t^{(n)}\circ \pi^{(n)}) -f_t^{(n)}|^2}\, \d t=0.$$ Combining this with \eqref{4.11} and Proposition \ref{P3.2}, we obtain
\beq\label{4.12}\lim_{n\to\infty} \lim_{\vv_0\downarrow 0} \Big(|u_{s,\vv_0}^{\ll,n}-u_s^\ll| + \|\nn^{(2)}(u_{s,\vv_0}^{\ll, n} - u_s^\ll)\| \Big)=0.\end{equation} Thus, taking inner product for both sides of \eqref{4.9} with every $e_i$ and letting $\vv_0\downarrow 0$ and $n\uparrow\infty$, we conclude that
  $\P$-a.s. for all $t\in [0, T\land\tau]$ and all $i\ge 1$,
\beg{equation*}\beg{split}
& \<u_t^\ll(Z_t), e_i\> -\e^{-\ll_i t} \<u_0^\ll(Z_0),e_i\> =
 \int_0^t \e^{-\ll_i(t-s)}\ll_i\<u_s^\ll(Z_s),e_i\>\d s \\
 &+ \int_0^t \e^{-\ll_i(t-s)} \big\<(\nn_{\si_s\d W_s}^{(2)}u_s^\ll)(Z_s), e_i\big\>
 +\int_0^t \e^{-\ll_i(t-s)} \<\ll u_s^\ll(Z_s)-b_s(Z_s),e_i\>\d s.
\end{split}
\end{equation*}
That is, $\P$-a.s. for all $t\in [0, T\land\tau]$,
\beg{equation*}\beg{split}
& u_t^\ll(Z_t) -\e^{tA_2}u_0^\ll(Z_0)=
 \int_0^t A_2\e^{(t-s)A_2}u_s^\ll(Z_s)\d s \\
 &+ \int_0^t \e^{(t-s)A_2}(\nn_{\si_s\d W_s}^{(2)}u_s^\ll)(Z_s)
 +\int_0^t \e^{(t-s)A_2}(\ll u_s^\ll(Z_s)-b_s(Z_s))\d s,
\end{split}
\end{equation*}
which then gives \eqref{4.2} by combining the second equation in \eqref{4.1}.
\end{proof}

As a consequence of Theorem \ref{T4.1}, we have the following  pathwise uniqueness result.

\beg{cor}\label{C4.2} Assume {\bf (H1)}-{\bf (H4)} and let $b: [0,\infty) \to  C_b(\H;\H_2)$ be measurable   such that
$$
\sup_{t\in [0,T]} \|b_t\|_{\gamma_\aa,\phi}<\infty,\ \ T>0
$$
holds for some $\aa\in (\ff 2 3,1)$ and $\phi\in\D_1$.  Let $(X_t,Y_t)_{t\ge 0}$ and $(\tt X_t,\tt Y_t)_{t\ge 0}$ be two adapted continuous processes on $\H$   with $(X_0,Y_0)=(\tt X_0,\tt Y_0).$ For any $n\ge 1$ let
$$\tau_n =n\land \inf\{t\ge 0: |(X_t,Y_t)|\ge n\},\ \ \  \tt \tau_n =n\land \inf\{t\ge 0: |(\tt X_t, \tt Y_t)|\ge n\}.$$ If $(X_t,Y_t)$ and $(\tt X_t,\tt Y_t)$ are mild solutions to $\eqref{1.1}$ for $t\in [0,\tau_n\land \tt \tau]$,   then  $\tau_n=\tt \tau_n$ and $(X_t,Y_t)=(\tt X_t,\tt Y_t)$ for all $t\in [0, \tau_n].$

When $\H_1,\H_2$ are finite-dimensional, the assertion holds for $\D_2$ in place of $\D_1$.   \end{cor}

\beg{proof} Due to Theorem \ref{T4.1} and Proposition \ref{P3.1}, the proof is similar to that of \cite[Proposition 3.1]{W14}. We address it here for complement.
Write
$$
T_n:=\tau_n\land\tt\tau_n, \ Z_t:= (X_t,Y_t),\  \tt Z_t:= (\tt X_t,\tt Y_t).
$$
It suffices to prove that for any $T>0$,
 \beq\label{A1} \int_0^T\E\big[1_{\{s<T_n\}}|Z_{s}-\tt Z_{s}|^2 \big]\d s =0.\end{equation}    By Theorem \ref{T4.1} for $\tau= T_n$, we have
 \beg{equation*} \beg{split} h(t)&:=\E\big[1_{\{t<T_n\}} |Z_{t}-\tt Z_{t}|^2\big]= \E \big[1_{\{t<T_n\}}\big(|X_{t}-\tt X_{t}|^2+ |Y_{t}-\tt Y_t |^2\big)\big]\\
 &\le C_1\bigg[\E\int_0^{t} 1_{\{s<T_n\}}|Y_s-\tt Y_s|^2\d s +\lambda\E\int_0^{t}1_{\{s<T_n\}} |u_s^\ll( Z_s)-u_s^\ll(\tt Z_s)|^2\d s \\
&\qquad\quad+I(t)+J(t)\bigg]+ 2\E\big[1_{\{t<T_n\}}|u_{t}^\ll(Z_{t})- u_{t\land T_n}^\ll(\tt Z_{t\land T_n})|^2\big]
 \end{split}\end{equation*} for some constant $C_1>0$ and $t\in [0,T]$,
 where
 \beg{equation}\label{IJ} \beg{split} &I(t):= \E \bigg|1_{\{t<T_n\}}\int_0^{t}A_2\e^{(t-s)A_2}\big(u_s^\ll(Z_s)-u_s^\ll(\tt Z_s)\big)\d s\bigg|^2,\\
 & J(t):= \E \int_0^{t} 1_{\{s<T_n\}}\|\e^{(t-s)A_2}\|_{HS}^2 \big\|\nn^{(2)} u_s^\ll(Z_s)-\nn^{(2)} u_s^\ll(\tt Z_s)\big\|^2\d s.\end{split}\end{equation}
 B \eqref{3.2'}, we have
 \beq\label{A}
|u^\ll_{s}(Z_{s})-u_{s}^\ll(\tt Z_{s})|
+\|\nn^{(2)} u^\ll_{s}(Z_{s})-\nn^{(2)} u_{s}^\ll(\tt Z_{s})\|\le \theta(\lambda)|Z_{s}-\tt Z_{s}|.
 \end{equation}
This implies that for all $r\in [0,T]$,
\beq\label{A2}
g(r):= \int_0^r h(t)\d t\le C_2\int_0^r g(t)\d t +C_2\int_0^r\big(I(t)+J(t)\big)\d t+2\theta(\lambda)^2g(r).
\end{equation}
By \eqref{IJ} and H\"older's inequality, we have
 \beg{equation*}\beg{split}I(t) &=  \sum_{i\ge 1} \E \bigg|1_{\{t<T_n\}}\int_0^{t}\ll_i\e^{-(t-s)\ll_i} \big\<u_s^\ll(Z_s)-u_s^\ll(\tt Z_s), e_i\big\> \d s\bigg|^2\\
 &\le \sum_{i\ge 1}  \E \int_0^{t}1_{\{s<T_n\}}\ll_i \e^{-(t-s)\ll_i}\<u_s^\ll(Z_s)-u_s^\ll(\tt Z_s), e_i\>^2 \d s,
 \end{split}\end{equation*}
which, together with Fubini's theorem and  \eqref{A}, yields
 \beq\label{A4} \beg{split}
\E\int_0^r I(t)\d t
 &\le \sum_{i\ge 1}\ll_i \E\int_0^r\d t  \int_0^{t}  1_{\{s<T_n\}} \e^{-(t-s)\ll_i}\<u_s^\ll(Z_s)-u_s^\ll(\tt Z_s), e_i\>^2 \d s\\
 &= \sum_{i=1}^\infty \E\int_0^{r}1_{\{s<T_n\}}\<u_s^\ll(Z_s)-u_s^\ll(\tt Z_s), e_i\>^2 \d s \int_s^{r}\ll_i \e^{-(t-s)\ll_i}\d t\\
 &\le  \E\int_0^r1_{\{s<T_n\}} |u_{s}^\ll(Z_{s})- u_{s}^\ll(\tt Z_{s})|^2 \d s \\
 &\le \theta(\lambda)^2 \int_0^r\E\big[1_{\{s<T_n\}}|Z_{s}-\tt Z_{s}|^2
\big]\d s
 =\theta(\lambda)^2 g(r).
 \end{split}\end{equation}
On the other hand, by \eqref{IJ}, \eqref{A} and \eqref{EW1}, we have
\beg{equation*}\beg{split}
\E\int_0^r J(t)\d t
&= \E \int_0^{r} 1_{\{s<T_n\}}\|\nn^{(2)} u^\ll_{s}(Z_{s})-\nn^{(2)} u^\ll_{s}(\tt Z_{s})\|^2 \d s \int_{s}^r\|\e^{(t-s)A_2}\|_{HS}^2\d t \\
&\le c_2\theta(\lambda)^2 \int_0^r h(s) \d s= c_2\theta(\lambda)^2 g(r).
\end{split}\end{equation*}
Substituting this and   \eqref{A4} into \eqref{A2}, we arrive at
 $$
 g(r)\le C_2 \int_0^r g(t)\d t+(C_2+2+c_2)\theta(\lambda)^2 g(r),\ \ \ r\in [0,T].
 $$
Therefore, by letting $\lambda$ be large enough so that $(C_2+2+c_2)\theta(\lambda)^2<\frac{1}{2}$,
then using Gronwall's inequality, we obtain  $g(T)=0$. The proof of the first assertion is complete.

Finally, when $\H_1$ and $\H_2$ are finite-dimensional and the condition on $b$ holds for $\phi\in \D_2$,
in the above proof, we only need to give a different treatment for $\E\int_0^r J(t)\d t$. Notice that in this case, $\sup_{t\ge 0}\|\e^{tA_2}\|_{HS}<\infty$,
and by taking $r= 1\land \ss{g(t)}$ in \eqref{3.2}, we obtain
\beq\label{BBC} \beg{split}&\|\nn^{(2)} u_{s}(Z_{s})-\nn^{(2)} u_{s}(\tt Z_{s})\|\\
&\le C_3 \bigg\{\ss{g(t)} +|Z_{s}-\tt Z_{s}| \bigg(1+\int_{g(t)^{\delta/2}\land 1}^1 \frac{\phi(s')}{s'}\d s'\bigg)\bigg\},\ s\in [0,T], t> s
\end{split}\end{equation}
for  some constant $C_3>0$. Hence, for some constants $C_4,C_5>0$, we have
\beg{equation*}\beg{split}
\E\int_0^r J(t)\d t
&= \E \int_0^{r}\d t \int_0^{t} 1_{\{s<T_n\}}\|\e^{(t-s)A_2}\|_{HS}^2 \big\|\nn^{(2)} u_s^\ll(Z_s)-\nn^{(2)} u_s^\ll(\tt Z_s)\big\|^2\d s \\
&\le C_4  \int_0^{r} \d t\int_0^t \bigg\{g(t)+ h(s)  \bigg(1+\int_{g(t)^{\delta/2}\land 1}^1 \frac{\phi(s')}{s'}\d s'\bigg)^2\bigg\} \d s\\
&\le C_5 \int_0^r g(t) \bigg(1+\int_{g(t)^{\delta/2}\land 1}^1 \ff{\phi(s)}s\d s\bigg)^2\d t,\ \ \ r\in [0,T].
\end{split}\end{equation*}
Substituting this and   \eqref{A4} into \eqref{A2} for large $\ll$, we arrive at
 $$
 g(r)\le C_6 \int_0^r g(t) \bigg(1+\int_{g(t)^{\delta/2}\land 1}^1 \ff{\phi(s)}s\d s\bigg)^2\d t=\int^r_0\rho(g(t))\d t,\ \ \ r\in [0,T]
 $$
 for some constant $C_6>0$, where $\rho(t):=C_6t\left(1+\int_{t^{\delta/2}\land 1}^1 \ff{\phi(s)}s\d s\right)^2$.
 Since $\phi\in \D_2$, we have $$\int_0^1 \ff{\d t}{\rho(t)} =\frac{2 }{C_6\delta} \int_0^1 \ff{\d t}{t \big(1+\int_{t\land 1}^1 \ff{\phi(s)}s\d s\big)^2}=\infty.$$ Therefore, by Bihari's inequality, this implies  $g(T)=0$  which is equivalent to  \eqref{A1}.
 \end{proof}

\section{Proofs of Theorem \ref{T1.1} and Theorem \ref{T1.2}}

Due to Corollary \ref{C4.2}, the proofs are more or less similar to \cite[\S 5]{W14}.

\beg{proof}[Proof of Theorem \ref{T1.1}] We split the proof into the following three steps.

(a) Firstly, we assume that for any $T>0$ there exist $\aa\in (\ff 2 3,1)$ and $\phi\in \D_1$   such that
$\sup_{t\in [0,T]}\|b_t\|_{\gamma_\aa,\phi}<\infty.$
By Girsanov's theorem, we see that for any $T>0$ and any initial value $(X_0,Y_0)$,
the equation \eqref{1.1} has a weak mild solution up to time $T$. Moreover, by Corollary \ref{C4.2},
the mild solution is pathwise unique up to arbitrary time $T>0$. So,  by the Yamada-Watanabe principle \cite{YW}
  (see \cite[Theorem 2]{OD04} or \cite{LR13} for the result in infinite dimensions), for any initial value,
  the equation \eqref{1.1}  has a unique mild solution which is non-explosive (i.e. the life time $\zeta=\infty$).

(b) In general,  take $\psi\in C_b^\infty ([0,\infty))$ such that $0\le \psi\le 1, \psi(r)=1$ for $r\in [0,1]$ and $\psi(r)=0$ for $r\ge 2.$  For any $m\ge 1,$  let
$$  b_t^{[m]}(z)= b_{t\land m}(z) \psi(|z|/m),\ \  t\ge 0,z\in\H.$$ Then by the condition in Theorem \ref{T1.1}(1), for any $m\ge 1$ and $T>0$, there exist
$\aa\in (\ff 2 3,1)$ and $\phi\in \D_1$ such that
$$
\sup_{t\in [0,T]}\|b_t^{(m)}\|_{\gamma_\aa,\phi}<\infty.
$$
So, for any initial value $(X_0,Y_0)$, the equation \eqref{1.1} for $b^{[m]}$ in place of $b$ has a unique mild solution  $(X_t^{(m)},Y_t^{(m)})_{t\ge 0}$ which is non-explosive. Let
$$
\tau_m=m\land \inf\big\{t\ge 0: |(X_t^{(m)},Y_t^{(m)})|\ge m\big\},\ \ \ m\ge 1.
$$
Since   $b_s^{[m]}(z)=b_s(z)$  for $s\le m$ and $|z|\le m,$   by Corollary \ref{C4.2}, for any $n,m\ge 1$ we have
$(X_t^{(n)},Y_t^{(n)})=(X_t^{(m)},Y_t^{(m)})$ for $t\in [0, \tau_n\land\tau_m]$. In particular, $\tau_m$ is increasing in $m$.  Let $\zeta=\lim_{m\to\infty} \tau_m$ and
$$(X_t,Y_t)= \sum_{m=1}^\infty 1_{[\tau_{m-1}, \tau_m)}(t) (X_t^{(m)},Y_t^{(m)})\ \ \ \tau_0:=0, t\in [0,\zeta).$$ Then it is easy to see that $(X_t,Y_t)_{t\in [0,\zeta)}$ is a mild solution to \eqref{1.1} with life time $\zeta$  and,  due to Corollary \ref{C4.2},  the mild solution is unique.  Thus, the assertion (1) is proved.

(c) Let  \eqref{C1} hold for some positive increasing $\ell, h$ such that
$\int_1^\infty \ff{\d s}{\ell_t(s)} =\infty, t\ge 0.$ Let $(X_t,Y_t)_{t\in [0,\zeta)}$ be a mild solution to \eqref{1.1}  with life time $\zeta$. Let
$\xi_t=\int_0^t \e^{(t-s)A_2} \si_s\d W_s$. Due to the boundedness of $\si$ and {\bf (H3)}, $\xi_t$ is an adapted continuous process on $\H_2$ up to the life time $\zeta$ with
$$\eta_T= |Y_0|^2 + 2\int_0^{T\land\zeta} h_{T\land\zeta}(|\xi_s|)\d s<\infty, \ \ T> 0.$$  Then $\tt Y_t:= Y_t -\xi_t$ solves  the equation
$$\d \tt Y_t = \big(A_2 \tt Y_t + b_t (X_t, \tt Y_t+\xi_t)\big)\d t,\ \ \ \tt Y_0=\tt Y_0,\ \ t<\zeta.$$ Due to \eqref{C1}, the increasing property of $h,\ell$,  and   $A_2\le 0$, this implies that for any $T>0$,
\beg{equation*}\beg{split}
\d |\tt Y_t|^2 &\le 2\<\tt Y_t,  b_t(X_t, \tt Y_t+\xi_t)\>\,\d t\\
&\le 2\big(\ell_{T\land \zeta}(|X_t|^2+|\tt Y_t|^2)+ h_T(|\xi_t|)\big)\d t,\ \ |\tt Y_0|^2=|Y_0|^2, t<\zeta\land T.\end{split}\end{equation*} So,
$$|\tt Y_t|^2 \le \eta_T +2\int_0^t \ell_T(|X_r|^2+ |\tt Y_r|^2)\d r,\ \ \ T>0, t\in [0,T\land\zeta).$$
On the other hand, there exists a random variable $C>1$ such that
\beq\label{XX}
|X_r|^2 =\bigg|\e^{rA_1} X_0+ \int_0^r \e^{(r-s)A_1}B\tilde Y_s\d s\bigg|^2\le  C   + (C-1) \sup_{s\in [0,r]} |\tt Y_s|^2,\ \ r\in [0, T\land\zeta).
\end{equation}
Therefore,  $g(t):=\sup_{s\in [0,t]}|\tt Y_s|^2$ satisfies
$$g(t)\le \eta_T +2\int_0^t \ell_T(C + C g(r))\d r,\ \ \ T>0, t\in [0,T\land\zeta).$$ Letting
$$\Gamma_T(s)=\int_1^s\ff{\d r}{2\ell_{T\land \zeta}(C+Cr)},$$
by Biharis's inequality, we obtain
\beq\label{XY} g(t) \le \Gamma_{T}^{-1} \big(\Gamma_{T}(\eta_T)+t\big),\ \ T>0,   t\in [0,\zeta\land T).\end{equation}
This implies  $\P(\zeta<\infty)=0$, i.e. $(X_t,Y_t)$ is non-explosive. Indeed, by the property  of the life time and \eqref{XX} and \eqref{XY}, on the set $\{\zeta\le T\}$ we have $\P$-a.s.
   $$\infty= \limsup_{t\uparrow\zeta}(|X_t|^2+ |Y_t|^2 ) \le \limsup_{t\uparrow\zeta}g(t)\le  \Gamma_{T}^{-1} \big(\Gamma_{T}(\eta_T)+T\big)<\infty,$$ where the last step is due to the fact that $\Gamma_T (r)\uparrow\infty$ as $r\uparrow\infty$, which implies  $\Gamma_T^{-1}(r)<\infty$ for any $r\in (0,\infty).$  This    contradiction   means that  $\P(\zeta\le T)=0$ holds for all $T\in (0,\infty)$. Hence, $\P(\zeta <\infty)=0.$
\end{proof}

\beg{proof}[Proof of Theorem \ref{T1.2}] Since $\H$ is finite-dimensional, we may apply It\^o's formula for $|X_t|^2$ directly to prove the non-explosion from \eqref{C1'} as in (c) in the proof of Theorem \ref{T1.1}.  So,  as explained in the proof of Theorem \ref{T1.1}, it remains to prove the pathwise uniqueness  for  $b$ satisfying
$$\sup_{t\in [0,T]} \|b_t\|_{\gamma_\aa,\phi}<\infty, \ \ T>0,$$ where $\aa\in (\ff 2 3, 1)$ and $\phi\in \D_2$. To apply Corollary \ref{C2.4}, we reformulate \eqref{1.2} as \eqref{1.1} for $A_1=A, A_2= I_d, $ and $b_t-I_d$ in place of $b_t$, where $I_d$ is the identity operator on $\R^d$. Then assumptions
{\bf(H1)}, {\bf (H3)} and {\bf (H4)} are trivial since $\H$ is finite-dimensional. Moreover, {\bf (H2)} holds for $A_0= I_m-A$. Thus, the pathwise uniqueness follows from Corollary \ref{C2.4}.\end{proof}

%\paragraph{\bf Acknowledgement.} The author would like to thank    for helpful comments and corrections.

\beg{thebibliography}{99}

\bibitem{DF}  G. Da Prato, F. Flandoli (2010), Pathwise uniqueness for  a class of SDE in Hilbert spaces and applications, \emph{J. Funct. Anal.} 259: 243--267.

\bibitem{DR1}   G. Da Prato, F. Flandoli, E. Priola, M. R\"ockner (2013), Strong uniqueness for stochastic evolution equations in Hilbert spaces perturbed by a bounded measurable drift, 
  \emph{Ann. Probab.} 41: 3306--3344.

\bibitem{DR2}  G. Da Prato, F. Flandoli, E. Priola, M. R\"ockner,  Strong uniqueness for stochastic evolution equations with  unbounded measurable drift term, to appear in \emph{J. Theor. Probab.} DOI: 10.1007/s10959-014-0545-0.

\bibitem{RN}   G. Da Prato, F. Flandoli, M. R\"ockner,  A. Yu. Veretennikov,  Strong uniqueness for SDEs in Hilbert spaces with non-regular drift, 
  arXiv: 1404.5418.

\bibitem{DZ}  G. Da Prato, J. Zabczyk (1992),  Stochastic Equations in Infinite Dimensions, \emph{Cambridge University Press}, Cambridge.

%\bibitem{FP} F. Flandoli, M. Gubinelli, E. Priola, \emph{Well-posedness of transport equation by stochastic perturbation,} Invent. Math. 180(2010), 1--53.

%\bibitem{FGP}  F. Flandoli, M. Gubinelli, E. Priola, \emph{Flow of diffeomorphisms for SDEs with unbounded H\"older continuous drift,} Bull. Sci. Math. 134(2010), 405--422.

%\bibitem{GP} I. Gy\"ongy, E. Pardoux, \emph{On the regularization effect of space-time white noise on quasi-linear parabolic partial differential equations,} Probab. Theory Relat. Fields 97(1993), 211--229.

\bibitem{GW}A. Guillin, F.-Y. Wang (2012),  Degenerate Fokker¨CPlanck equations: Bismut formula, gradient estimate and Harnack inequality, \emph{J. Differential
Equations} 253: 20¨C40.

\bibitem{KR} N. V. Krylov, M. R\"ockner (2005), Strong solutions of stochastic equations with singular time dependent drift, \emph{Probab. Theory Relat. Fields}  131: 154--196.

\bibitem{OD04} M. Ondrej\'et (2004),  Uniqueness for stochastic evolution equations in Banach spaces, \emph{Dissertationes Math. (Rozprawy Mat.) }

\bibitem{Lu} A. Lunardi (1995), Analytic Semigroups and Optimal Regularity in Parabolic Problems, \emph{Birkh\"auser}, 1995, Switzerland.

\bibitem{LR13}  C. Pr\'ev\^ot,  M. R\"ockner (2007), A Concise Course on Stochastic   Partial Differential Equations, \emph{Lecture Notes in Math. } Vol. 1905, Springer, Berlin, 2007.

%\bibitem{RW10} 	M. R\"ockner, F.-Y. Wang, Log-Harnack  Inequality for Stochastic differential equations in Hilbert spaces and its consequences, Infinite Dimensional Analysis, Quantum Probability and Related Topics 13(2010), 27--37.

%\bibitem{Wbook} F.-Y. Wang,  \emph{Harnack Inequalities for Stochastic Partial Differential Equations,} Springer, Berlin, 2013.

%\bibitem{Wb2} F.-Y. Wang,  \emph{Analysis for Diffusion Processes on Riemannian Manifolds,} World Scientific, Singapore, 2013.

\bibitem{W14}   F.-Y. Wang,  Gradient Estimates and Applications for SDEs in Hilbert Space with Multiplicative
Noise and Dini Drift,  http://arxiv.org/pdf/1404.2990v3.pdf.

%\bibitem{WZ14}  	F.-Y. Wang, T. Zhang, \emph{Log-Harnack inequalities for semi-linear SPDE with strongly multiplicative noise,} Stoch. Proc. Appl. 124(2014), 1261--1274.

\bibitem{YW} T. Yamada, S. Watanabe (1971), On the uniqueness of solutions of stochastic differential equations, \emph{J. Math. Kyoto Univ.} 11: 155--167.

\bibitem{V} A. J. Veretennikov (1980),  Strong solutions and explicit formulas for solutions of stochastic integral equations,  \emph{(Russian) Mat. Sb. (N.S.)} 111: 434--452.

\bibitem{Zh} X. Zhang (2005), Strong solutions of SDES with singular drift and Sobolev diffusion coefficients, \emph{Stoch. Proc. Appl. } 115: 1805--1818.

\bibitem{Z} A. K. Zvonkin (1974), A transformation of the phase space of a diffusion process that will remove the drift, \emph{(Russian) Mat. Sb. (N.S.)} 93: 129--149.
\end{thebibliography}
\end{document}